\documentclass[11pt,a4paper,fleqn]{article}
\usepackage{amsfonts,amsmath}
\usepackage{latexsym}
\usepackage{amssymb}
\usepackage{euscript}
\usepackage{graphicx}
\usepackage{a4wide}
\usepackage{epsfig}

\newtheorem{prop}{Proposition}[section]

\newtheorem{lemma}[prop]{Lemma}

\newtheorem{rem}[prop]{Remark}
\newtheorem{theo}[prop]{Theorem}

\newcommand{\noi}{\noindent}
\newcommand{\dis}{\displaystyle }

\newenvironment{proof}[1]{\begin{trivlist}\item {\it
\bf Proof.}\quad} {\qed\end{trivlist}}
\newenvironment{prooff}[1]{\begin{trivlist}\item {\it
\bf Proof}\quad} {\qed\end{trivlist}}

\newcommand{\qed}{\nopagebreak\hspace*{\fill}
{\vrule width6pt height6ptdepth0pt}\par}

\begin{document}


\title{\bf LIMITING LAWS FOR LONG  BROWNIAN BRIDGES PERTURBED
BY THEIR ONE-SIDED MAXIMUM,  III}
\author{{\small{\text{\bf Bernard ROYNETTE}$^{(1)}$}},\
{\small {\text{\bf  Pierre VALLOIS}$^{(1)}$}} {\small and} {\small
{\text{\bf Marc YOR }$^{(2),(3)}$}}}



\maketitle {\small

\noindent (1)\,\, Universit\'e Henri Poincar\'e, Institut de
Math\'ematiques Elie Cartan, B.P. 239, F-54506 Vand\oe uvre-l\`es-Nancy Cedex\\

\noi (2)\,\, Laboratoire de Probabilit\'{e}s et
 Mod\`{e}les Al\'{e}atoires, Universit\'{e}s Paris VI et VII -  4, Place Jussieu
 - Case 188 -
 F-75252 Paris Cedex 05.\\

 \noi (3) Institut Universitaire de France.\\


\noi \begin{flushright}In homage to Professors E. Csaki and P.
Revesz.
\end{flushright}

 \vskip 40 pt \noi {\bf Abstract.} Results of penalization of a
 one-dimensional Brownian motion $(X_t) $, by its one-sided
 maximum $\dis (S_t=\sup_{0 \leq u \leq t}X_u)$, which were
 recently obtained by the authors are improved  with the consideration-in the
 present paper- of the asymptotic behaviour of the likewise
 penalized Brownian bridges of length $t$, as $t\rightarrow
 \infty$, or penalizations by functions of $(S_t,X_t)$, and also the
 study of the speed of convergence, as $t\rightarrow
 \infty$, of the penalized distributions at time $t$.


\vskip 40pt


\noi {\bf Key words and phrases} :  penalization, one-sided
maximum, long Brownian bridges, local time, Pitman's theorem

\smallskip



\noi {\bf AMS 2000 subject classifications} :  60 B 10, 60 G 17,
60 G 40, 60 G 44, 60 J 25, 60 J 35, 60 J 55, 60 J 60, 60 J 65.

\section{Introduction}\label{int}

\setcounter{equation}{0}
 \noi {\bf 1.1} Let $\big(\Omega = {\cal C}(\mathbb{R} _+,
\mathbb{R}),$ $ (X_t)_{t \geq 0},\; ({\cal F}_t)_{t \ge 0}\big)$
be the canonical space with $(X_t)$ the process of coordinates :
$X_t(\omega)=\omega(t) ; t\geq 0$, $({\cal F}_t)_{t \ge 0}$ the
canonical filtration associated with $(X_t)$. We write  ${\cal
F}_\infty$ for the  $\sigma$-algebra generated by $\dis \bigcup_{t
\ge 0}{\cal F}_t$. Let $P_0$ be
  the  Wiener measure defined on the
canonical space such that $P_0(X_0=0)=1$.

\noi In this paper, as well as in the previous ones
(\cite{RoyValYor}, \cite{RVY2}, \cite{RVY2a}), we consider
perturbations of Brownian motion with certain processes $( F_t)_{t
\ge 0}$, which we call weight-processes; precisely, let $( F_t)_{t
\ge 0}$
   be an $({\cal F}_t)$-adapted,
non negative process, such that $0 < E_0 (F_t) < \infty$, for any
$t \ge 0$, and $Q_{0,t}^{F}$ the probability measure (p.m.)
defined on $(\Omega, \; {\cal F}_t)$ as follows :
\begin{equation}\label{int1}
    Q_{0,t}^{F} (\Gamma_t) := \frac{1 }{ E_0 [F_t]} \; E_0[1_{\Gamma_t} F_t],
\quad  \Gamma_t \in {\cal F}_t.
\end{equation}
%
We can interpret   the p.m. $Q_{0,t}^F$ as the Wiener measure
penalized by the weight $ F_t$. We say that a penalization
principle holds if there exists a p.m. $Q_0^F$ on $\big(\Omega,
{\cal F}_\infty\big)$ such that $Q_{0,t}^F$ converges weakly to
$Q_0^F$,  as $t\rightarrow \infty$ :
\begin{equation}\label{int2}
    \lim_{t\rightarrow \infty}Q_{0,t}^F(\Gamma_u)=Q_0^F
    (\Gamma_u),\quad \mbox{for any} \ \Gamma_u \in {\cal F}_u, \
    u\geq 0.
\end{equation}

\noi Throughout the paper, $(S_t)$ stands for the one-sided
maximum of $(X_t)$ : $\dis S_t:= \max_{0\leq u \le t}X_u, \ t\geq
0$.

\noi In fact, in our study, the following situation always occurs
: let
$$M_u^{(t)}:=\frac{1}{E_0[F_t]}E_0[F_t|{\cal F}_u], \quad u <t.$$
Then, we show that, for fixed $u$, $M_u^{(t)}$ converges a.s.,
with respect to $P_0$, to a variable $M_u$, such that
$E_0[M_u]=1$. Thus by Scheff\'{e}'s lemma (see, e.g. \cite{Mey},
Chap. V,  T21) $M_u^{(t)}$ converges in $L^1(P_0)$ towards $M_u$,
which explains why (\ref{int2}) holds without any restriction on
$\Gamma_u \in {\cal F}_u$.

 \noi {\bf 1.2} In a series of papers
(\cite{RoyValYor}, \cite{RVY2}, \cite{RVY2a} and  \cite{RVY3}) we
have considered some classes of examples involving respectively
for our weight-process $(F_t)$ a function of :

\begin{itemize}
\item $\dis \int_0^tV(X_s)ds$ where $V :
\mathbb{R}\mapsto\mathbb{R}_+$.
    \item the unilateral maximum  $S_t$; we have also treated  the two-dimensional process
$(S_t,t)$.
    \item $(L_t^0; t \ge 0)$ the local time at $0$ of $(X_t)_{t \ge
0}$.

\item The triple $((S_t,I_t,L^0_t); t\geq 0)$, where $(I_t)$
denotes the one-sided minimum :  $\displaystyle I_t = - \min_{0
\le u \le t} X_u$.

    \item $(D_t; t \ge 0)$ the number of down-crossings of $X$ from level
$b$ to level $a$.
\end{itemize}

\noi In this paper we only consider the case : $F_t=f(X_t,S_t)$,
where $f : \mathbb{R}\times\mathbb{R}_+\mapsto \mathbb{R}_+$. In
particular if  $F_t=\varphi(S_t)$, where $\varphi :
\mathbb{R}_+\mapsto \mathbb{R}_+$ defines a probability density,
i.e.
\begin{equation}\label{int3}
    \int_0^\infty \varphi (y)dy=1 ,
\end{equation}
\noi our starting point is the following main result in
\cite{RVY2a} :

\begin{theo} \label{tint1}Let $\varphi : \mathbb{R}_+\mapsto \mathbb{R}_+$ satisfying (\ref{int3})
 and $\dis \Phi(y)=\int_0^y\varphi (x)dx, y\geq 0$.
\begin{enumerate} \item  For every  $u \ge 0$, and
$\Gamma_u$ in ${\cal F}_u$, the quantity  :
\begin{equation}\label{int4}
   Q^\varphi _0( \Gamma_u):=\lim_{t \to \infty} \; \frac{E_0\big[1_{\Gamma_u} \varphi
(S_t) \big] }{ E_0 \big[\varphi (S_t)\big]} ,
\end{equation}
exists; hence, $Q^\varphi _0$ may be extended as a p.m. on
$\big(\Omega, {\cal F}_\infty\big)$.

 \item  It is equal to $E_0 [1_{\Gamma_u} M_u^{\varphi}]$,
  where $(M_u^{\varphi})_{u \ge 0}$ is the
martingale :
\begin{equation}\label{int5}
    M_u^{\varphi}=\varphi(S_u)(S_u-X_u)+1-\Phi(S_u); \quad
    u\geq 0.
\end{equation}
\noi (These  $\big(P_0,({\cal F}_u)\big)$-martingales have been
introduced in \cite{AY1}).

\item The probability $Q^\varphi_0$ may be disintegrated as
follows :

\begin{enumerate}
    \item under $Q^\varphi_0$, $S_\infty$ is  finite a.s., and
    admits $\varphi$ as a probability density;

   \item  $ Q_0^\varphi (S_\infty \in dy)$  a.e., conditionally on
$S_\infty=y$, the law of $(X_t)$, under $ Q_0^\varphi$ is equal to
$Q^{(y)}_0$, where, for any $y>0$,  the p.m. $Q_0^{(y)}$  on the
canonical space is defined as follows :

\begin{enumerate}
        \item $(X_t; t \le T_y)$ is a Brownian motion started at
$0$, and considered up to $T_y$, its first hitting time  of $y$,
        \item the process $(X_{ T_{y}+t} ; t \ge 0)$ is a
         "three dimensional Bessel process below $y$", namely :
          $(y-X_{ T_{y}+t} \;;\; t \ge 0)$ is a three
dimensional Bessel process started at $0$.
        \item the processes $(X_t; t \le T_y)$ and $(X_{ T_{y}+t} ; t \ge
        0)$ are independent.
    \end{enumerate}

    \item Consequently :
    \begin{equation}\label{int8}
    Q_0^\varphi(\Gamma|S_\infty=y):=Q_0^{(y)}(\Gamma), \quad \mbox{ for
    any} \  \Gamma \in {\cal F}_\infty,
\end{equation}
\begin{equation}\label{int6}
    Q_0^{\varphi}(\cdot ) = \int_0^{\infty} Q^{(y)}_0(\cdot ) \; \varphi
(y) \, d y.
\end{equation}
\end{enumerate}
\end{enumerate}

\end{theo}

\noi In the present paper, we develop a number of variants of this
Theorem \ref{tint1}, by presenting either extensions or some new
proofs of this theorem. Here are these variants, together with the
organization of our paper.

\noi In Section 2, we give, in particular, another proof of
Theorem \ref{tint1}, which originates from the following
considerations : the main step in \cite{RVY2a} consisted in
studying the asymptotics of $E[\varphi(S_t)|{\cal F}_s]$, for
fixed $s$, as $t\rightarrow \infty$. In Section 2 here, we proceed
 in a dual manner by studying the asymptotics of
\begin{equation}\label{int7}
Q_{0,t}^{(y)}(\Gamma_u):=P(\Gamma_u|S_t=y),
\end{equation}
as $t\rightarrow \infty$, where $u\geq 0$ and $\Gamma_u \in {\cal
F}_u$ are fixed.

\begin{theo}\label{tint2} Let $y>0$, $u\geq 0$ and $\Gamma_u \in {\cal
F}_u$.
\begin{enumerate}
    \item As $t\rightarrow \infty$, $Q_{0,t}^{(y)}(\Gamma_u)$
    converges towards the probability $Q^{(y)}_0(\Gamma_u)$, where
$Q^{(y)}_0$ is the probability introduced in Theorem \ref{tint1},
3.
    \item Moreover, $Q^{(y)}_0$ satisfies :
    \begin{equation}\label{int9}
    Q^{(y)}_0(\Gamma_u)=e^{-y^2/2u}\sqrt{\frac{2}{\pi
    u}}E_0\Big[1_{\Gamma_u}(y-X_u)\Big|S_u=y\Big] +
    E_0[1_{\Gamma_u}1_{\{S_u<y\}}].
\end{equation}
\end{enumerate}

\end{theo}

\noi In Section \ref{bri}, we strengthen the result obtained in
Section \ref{prt}, in that we consider the existence of the
limits, as $t\rightarrow \infty$, of :
\begin{equation}\label{int10}
    \frac{E_0\big[1_{\Gamma_u} \varphi
(S_t) |X_t=a\big] }{ E_0 \big[\varphi (S_t)\big|X_t=a]}
\end{equation}
and, in the spirit of the preceding Section \ref{prt} (or Theorem
\ref{tint2}) :
\begin{equation}\label{int11}
    Q_{0,t}^{a,y}(\Gamma_u):=P_0(\Gamma_u |X_t=a,S_t=y),
\end{equation}
where $u\geq 0,  \Gamma_u \in {\cal F}_u, y\geq a_+$.

\noi The title of the present paper originates from this  central
Section \ref{bri}. The results are the following :
\begin{itemize}
    \item concerning (\ref{int11}), we obtain :

\begin{theo}\label{tint3} \begin{enumerate}
    \item  For any $u\geq 0$ and $ \Gamma_u \in {\cal F}_u$,
    \begin{equation}\label{int12}
    \lim_{t\rightarrow
    \infty}Q_{0,t}^{a,y}(\Gamma_u):=Q_{0}^{a,y}(\Gamma_u),
\end{equation}
exists.
    \item The p.m. $Q_0^{a,y}$ may be expressed as a convex
    combination of the laws $Q_0^{(z)}, z\in \mathbb{R}_+$ :
    \begin{equation}\label{int13}
    (2y-a)Q_{0}^{a,y}(\cdot)=(y-a)Q_{0}^{(y)}(\cdot)+\int_0^y dz
    Q_{0}^{(z)}(\cdot).
\end{equation}
\end{enumerate}

\end{theo}

\begin{rem}\label{rbri1}\begin{enumerate}
    \item Recall that $Q^{(y)}_0(S_\infty =y)=1$. Since
     $Q_{0}^{a,y}$  satisfies  (\ref{int13}), we deduce :
$$ Q_{0}^{a,y}(S_\infty =y)=\frac{y-a}{2y-a}.$$
    \item As we started with the Brownian  bridge, we might have expected that, under the
     limiting
    p.m., some constraint involving
     the position of the process
     at infinity would hold. This is not the case, indeed,
    the parameter $a$ only appears in the coefficients of the
    convex combination in (\ref{int13}) and $\dis Q^{(y)}_0\big( \lim_{t\rightarrow
\infty}X_t=-\infty \big)=1$.

\item Identity (\ref{int13}) implies that $(a,y)\mapsto
Q_{0}^{a,y}$ is continuous.

\item We may recover $Q_{0}^{(y)}$ from $\big(Q_{0}^{a,y};y_+\leq
a\big)$ since $\dis
Q_{0}^{(y)}=\frac{d}{dy}\big(yQ_{0}^{y,y}\big)$.

\item Let $\mu^{a,y}$ the p.m. on $\mathbb{R}_+$ : $\dis
\mu^{a,y}(dz)=\frac{y-a}{2y-a}\delta_y(dz)+\frac{1}{2y-a}1_{[0,y]}(z)dz$.
The relation (\ref{int13}) admits the following probabilistic
interpretation : first,   $z$ is chosen at random following
$\mu^{a,y}$; secondly, the dynamics of $(X_t)$ is given by
$Q_{0}^{(z)}$.

\item From L\'{e}vy's theorem, under $P_0$,
$\big((S_t-X_t,S_t;t\geq 0\big)$ and $\big((|X_t|,L^0_t;t\geq
0\big)$ have the same distribution. Let $\mathbb{Q}_0^{(y)}$ be
the unique p.m. on $\big(\Omega, \sigma(|X_t|, t\geq 0)\big)$
satisfying :
    \begin{equation}\label{int13c}
    \mathbb{Q}_0^{(y)}(\Gamma_u)=e^{-y^2/2u}\sqrt{\frac{2}{\pi
    u}}E_0\Big[1_{\Gamma_u}|X_u|\ \Big|L^0_u=y\Big] +
    E_0[1_{\Gamma_u}1_{\{L^0_u<y\}}],
\end{equation}
for any $u\geq 0$ and $\Gamma_u \in \sigma (|X_t|,t\leq u)$.

\noi In a forthcoming paper \cite{RVY4} it is proved that  the
analog of (\ref{int12}) and (\ref{int13}) is   :
    \begin{equation}\label{int13b}
\lim_{t\rightarrow \infty}P_0(\Gamma_u \Big|\
|X_t|=a,L^0_t=y)=\frac{a}{a+y}\mathbb{Q}_0^{(y)}(\Gamma_u)+\frac{1}{a+y}
\int_0^y \mathbb{Q}_0^{(z)}(\Gamma_u)dz,
\end{equation}
with  $\Gamma_u$  any event in $  \sigma (|X_t|,t\leq u)$, and an
adequate extension of this result with $|X_s|$ being replaced by a
Bessel process with dimension $d<2$ is obtained.

\end{enumerate}

\end{rem}

    \item As for (\ref{int10}), we obtain :
\begin{theo}\label{tint4} Let $\varphi : \mathbb{R}_+\mapsto
\mathbb{R}_+$ such that :
\begin{equation}\label{int14}
    \int_0^\infty (1+x)\varphi (x)dx < \infty.
\end{equation}
\begin{enumerate}\item For any $u\geq 0$, $ \Gamma_u \in {\cal F}_u$ and $a\in
\mathbb{R}$, we have :
\begin{equation}\label{int15}
  Q^{a,\varphi }_{0}( \Gamma_u):=
  \lim_{t \to \infty} \; \frac{E_0\big[1_{\Gamma_u} \varphi
(S_t) |X_t=a\big] }{ E_0 \big[\varphi (S_t)|X_t=a\big]} ,
\end{equation}
 exists.
\item The p.m. $Q^{a,\varphi }_{0}$ may be expressed in terms of
either of the two families $(Q^{a,y}_{0},y>0)$ and

\noi $(Q^{(y)}_{0},y>0)$ :
\begin{eqnarray}
Q^{a,\varphi }_{0}(\cdot)&=&\frac{1}{\int_{a_+}^\infty
(2y-a)\varphi(y)dy} \int_{a_+}^\infty (2y-a)\varphi(y) Q^{a,y}_{0}(\cdot)dy \label{int16} \\
&=& \dis \frac{1}{\int_{a_+}^\infty (2y-a)\varphi(y)dy}\Big[
\int_{a_+}^\infty (y-a)\varphi(y) Q^{(y)}_{0}(\cdot)dy \nonumber
\\
&& \qquad \qquad \qquad \qquad +\int_0^\infty \big(1-\Phi
\big(z\vee(a_+)\big)\big) Q^{(z)}_{0}(\cdot)dz \Big]
.\label{int17}
\end{eqnarray}
\end{enumerate}

\end{theo}

\end{itemize}

\noi We would like to generalize Theorem \ref{tint4}, by replacing
the weight-process $\big(\varphi (S_t)\big)$ with $\big(f(X_t
,S_t)\big)$, where $f : \mathbb{R}\times \mathbb{R}_+ \mapsto
\mathbb{R}_+$ is Borel.

\begin{theo}\label{tint5} To $f : \mathbb{R}\times \mathbb{R}_+ \mapsto
\mathbb{R}_+$ such that :
\begin{equation}\label{int18}
    \overline{f}:=\int_\mathbb{R} da \int_{a_+}^\infty (2y-a) f(a,y)dy < \infty
\end{equation}
we associate $f^\star=1/ \overline{f}$, and :
    \begin{equation}\label{int20}
    \varphi (y)=
    f^\star\Big[\int_\mathbb{R}da\int_{y\vee
    a_+}^\infty f(a,\eta)d\eta +\int _{-\infty}^yf(a,y)(y-a)da\Big].
\end{equation}

\begin{enumerate}
    \item For every  $u \ge 0$, and $\Gamma_u$ in ${\cal F}_u$,
    \begin{equation}\label{int19}
   \lim_{t \to \infty} \; \frac{E_0\big[1_{\Gamma_u} f(X_t,
S_t) \big] }{ E_0 \big[f(X_t,S_t)\big]}=Q_0^\varphi(\Gamma_u) ,
\end{equation}
where $Q_0^\varphi$ is the p.m. introduced in Theorem \ref{tint1},
associated with the $\big(P_0,({\cal F}_t)\big)$ martingale
$(M_t^{\varphi})$.
    \item Moreover the following relations hold :
\begin{equation}\label{int21}
    M_t^{\varphi}=f^\star \int_\mathbb{R}da\int_{
    a_+}^\infty (2y-a)f\big(a+X_t,S_t\vee(y+X_t)\big)dy,
\end{equation}
\begin{equation}\label{int22}
    Q_0^\varphi(\cdot)=f^\star\int_\mathbb{R}da\int_{
    a_+}^\infty (2y-a)f(a,y)Q^{a,y}_{0}(\cdot) dy,
\end{equation}
where the p.m. $Q_0^{a,y}$ is defined in Theorem \ref{tint3}.
\end{enumerate}

\end{theo}

\noi Theorem \ref{tint5} led us to go further and to enquire what
happens if $f : \mathbb{R}\times \mathbb{R}_+ \mapsto
\mathbb{R}_+$  does not satisfy (\ref{int18}). Rather than trying
 to give a complete answer, we shall restrict ourselves to
functions $f$ of  exponential type :
\begin{equation}\label{int23}
    f(a,y)= e^{\lambda y+\mu a}, \ y\geq a_+ , \lambda, \mu  \in
    \mathbb{R}.
\end{equation}
It is easy to check (see Section \ref{pexp}) that, if $f$ is given
by (\ref{int23}), then : $\overline{f}< \infty$ iff $ \mu
    >0 $ and $\ \lambda + \mu <0$. Then in this case Theorem
    \ref{tint5} applies.

\noi We claim that for any $\lambda, \mu  \in \mathbb{R}$ a
penalization principle holds and we are able to describe the
limiting p.m. Before stating this result in Theorem \ref{tint6}
below, let us introduce the three disjoint sets :
\begin{equation}\label{int24}
    R_1=\big\{ (\lambda , \mu)\in \mathbb{R}\times \mathbb{R}; \
    \lambda + \mu <0 , \mu \geq 0 \big \},
\end{equation}
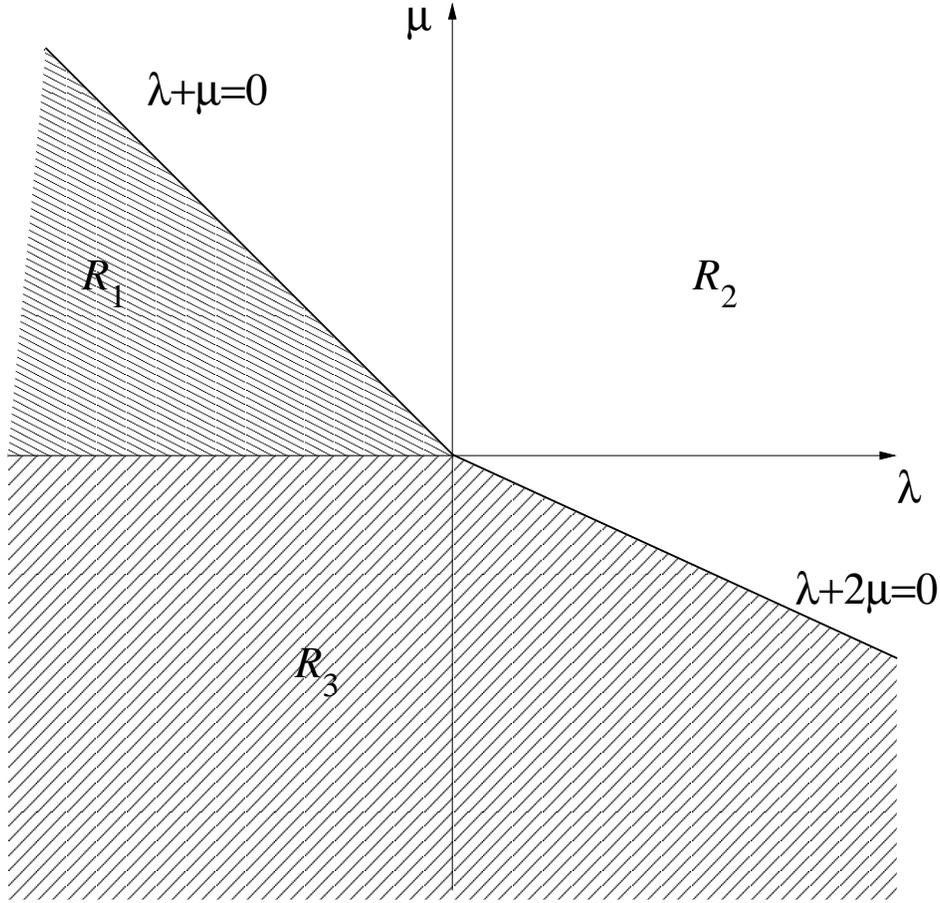
\begin{equation}\label{int25}
    R_2=\big\{ (\lambda , \mu)\in \mathbb{R}\times \mathbb{R}; \
    \lambda + 2 \mu \geq 0 ,\lambda + \mu \geq 0  \big \},
\end{equation}
\begin{equation}\label{int26}
    R_3=\big\{ (\lambda , \mu)\in \mathbb{R}\times \mathbb{R}; \
    \lambda + 2 \mu < 0 ,  \mu < 0 \big \}.
\end{equation}
%
See the figure below.

\begin{figure}[h]
\begin{center}
\includegraphics[width=0.8\textwidth]{pierre1.eps}
\end{center}
\end{figure}

\begin{theo} \label{tint6} Let $\lambda, \mu  \in \mathbb{R}$.
\begin{enumerate}
    \item For every  $u \ge 0$, and $\Gamma_u$ in ${\cal F}_u$,
    \begin{equation}\label{int27}
   \lim_{t \to \infty} \; \frac{E_0\big[1_{\Gamma_u} e^{\mu X_t +
\lambda S_t}) \big] }{ E_0 \big[e^{\mu X_t + \lambda S_t} \big]} ,
\end{equation}
exists and is equal to $E_0\big[1_{\Gamma_u}M_u^{\mu,
\lambda}\big]$, with $(M_u^{\mu, \lambda})$ a positive
$\big(({\cal F}_u),P_0\big)$ martingale, such that $M_0^{\mu,
\lambda}=1$, which is given by
\begin{equation}\label{int28}
M_u^{\mu, \lambda}= \left\{
\begin{array}{ll}
  -(\lambda + \mu)e^{(\lambda + \mu)S_u}(S_u-X_u)+e^{(\lambda + \mu)S_u} & \mbox{if}
  \ (\lambda , \mu)\in R_1, \\
  &\\
  e^{\{ (\lambda + \mu)X_u-(\lambda + \mu)^2u/2\}}& \mbox{if}
  \ (\lambda , \mu)\in R_2, \\
  &\\
 \dis  e^{\{ (\lambda + \mu)S_u- \mu^2u/2\}}\big[\cosh\big(\mu(S_u-X_u)\big)
  - \frac{\lambda + \mu}{\mu}\sinh\big(\mu(S_u-X_u)\big)\big]& \mbox{if}
  \ (\lambda , \mu)\in R_3. \\
\end{array}
\right.
\end{equation}
\item Consequently, $\Gamma_u (\in {\cal F}_u) \mapsto
E_0\big[1_{\Gamma_u}M_u^{\mu, \lambda }\big]$ induces a p.m. on
$\big(\Omega,{\cal F}_\infty\big)$.
\end{enumerate}

\end{theo}

\begin{rem} \label{rint1}\begin{enumerate}
    \item We have already observed that if $f$ is defined by (\ref{int23}),
     then $\overline{f}< \infty$ iff $ \mu
    >0 $ and $\ \lambda + \mu <0$. Thus, in this case, Theorem \ref{tint5} implies that  $(M_t^{\mu,
    \lambda})$ is a martingale of the type $(M^\varphi_t)$ where
    $\varphi$ is given by (\ref{int20}). An easy calculation yields :
     $\varphi (y)=-(\lambda + \mu) e^{(\lambda + \mu)y}, y\geq
    0$,  and :
    $$M_t^{\mu,\lambda}=M^\varphi_t=
    -(\lambda + \mu)e^{(\lambda + \mu)S_t}(S_t-X_t)+e^{(\lambda +
    \mu)S_t},\ t \geq 0.$$
    \item In the third case (i.e. $(\lambda , \mu)\in R_3$), the
    martingale belongs to the family of Kennedy martingales. These
    martingales were used in \cite{AY1} and play a central role in
    \cite{RVY2a}. Let us briefly recall the definition of these
    processes.

     \noi  To $\psi :\mathbb{R} \mapsto [0,\infty[$, a
    Borel function satisfying :
\begin{equation}\label{int29}
\int_{x}^\infty \psi(z)e^{-\lambda z}dz <\infty, \quad \forall
x\in \mathbb{R}.
\end{equation}
we associate  the function  $\Phi : \mathbb{R} \mapsto \mathbb{R}$
:
\begin{equation}\label{int30}
\Phi(y)=1-e^{\lambda y}\int _y^\infty \psi(z)e^{-\lambda
    z}dz, \ y \in\mathbb{R}.
\end{equation}
Let $\varphi$ be the derivative of $\Phi$; then, $\dis
\varphi(y):=\Phi'(y)=\psi(y)-\lambda e^{\lambda y}\int _y^\infty
\psi(z)e^{-\lambda
    z}dz$, and
\begin{equation}\label{int31}
  M_t^{\lambda,\varphi}:= \Big\{  \psi (S_t) \frac{\sinh \big(\lambda (S_t-X_t)\big)
}{ \lambda}+ e^{\lambda X_t}\int_{S_t}^\infty \psi (z)e^{-\lambda
z} dz\Big\}
  e^{- \lambda^{2}t/2},
\end{equation}
is a positive $\big(({\cal F}_t),P_0\big)$-martingale.

\item Let $Q_0^{\mu, \lambda}$ be the p.m. defined in point 2. of
Theorem \ref{tint6}, and $P_0^{\delta}$ be the law of Brownian
motion with drift $\delta$, starting at $0$. Using Theorem 3.9 of
\cite{RVY2a}, we may reformulate (\ref{int28}) as follows :

\begin{equation}
Q_0^{\mu, \lambda}= \left\{
\begin{array}{ll}
  Q_0^{\varphi_{-(\mu +\lambda)}} & \mbox{if}
  \ (\lambda , \mu)\in R_1, \\
  &\\
 P_0^{\mu+ \lambda}& \mbox{if}
  \ (\lambda , \mu)\in R_2, \\
  &\\
 \dis  \frac{\lambda + 2\mu}{2\mu}e^{\lambda S_\infty}\cdot P_0^{\mu}& \mbox{if}
  \ (\lambda , \mu)\in R_3 \\
\end{array}
\right.
\end{equation}
where $\varphi_{\delta}(y)=\delta e^{-\delta y}, \ \delta
>0,y\geq 0$.

\end{enumerate}

\end{rem}

\noi The proof of Theorem \ref{tint6} is postponed to Section
\ref{pexp}.

\noi Let $\varphi$ as in Theorem \ref{tint1}. We are now
interested in the rate of convergence of $\dis Q_{0,t}^\varphi
(\Gamma _u):=\frac{E_0[1_{\Gamma_u}\varphi (S_t)]}{E_0[\varphi
(S_t)]}$ towards $Q_{0}^\varphi (\Gamma _u)$, as $t\rightarrow
\infty$, for any $\Gamma_u\in {\cal F}_u$. More generally, under
additional assumptions, we are able to determine the asymptotic
development of $ Q_{0,t}^\varphi (\Gamma _u)$ in powers of $1/t$,
$t\rightarrow \infty$.

\begin{theo} \label{tint7}
Let $\varphi : \mathbb{R}_+\mapsto \mathbb{R}_+$ satisfying
(\ref{int3}) and the related function $\Phi$ as in Theorem
\ref{tint1}. We suppose that there exists an integer $n\geq 1$
such that :
\begin{equation}\label{int32}
    \int_0^\infty y^{2n+3}\varphi (y)dy<\infty.
\end{equation}
\begin{enumerate}
    \item There exists a family of functions $(F_i^\varphi)_{1\leq i \leq
    n}$, $F_i^\varphi  :  \mathbb{R}\times  \mathbb{R}_+\times
    \mathbb{R}_+ \mapsto \mathbb{R}$, such that
    \begin{enumerate}
        \item $(F_i^\varphi(X_t,S_t,t), t\geq 0)$ is a $\big(({\cal
        F}_t),P_0 \big)$-martingale, for any $1\leq i \leq n$,
        \item If $i=1$, we have :
        \begin{equation}\label{int33}
F_1^\varphi(X_t,S_t,t)=
-\widetilde{F}_1^\varphi(X_t,S_t)+\big(t+\int_0^\infty y^2 \varphi
(y)dy \big)M_t^\varphi,
\end{equation}
where
\begin{equation}\label{int34}
    \widetilde{F}_1^\varphi(a,y)=\varphi (y)\frac{(y-a)^3}{3!}+\frac{1}{2}\int_y^\infty \varphi
    (v)(v-a)^3dv, \quad t,y \geq 0, x\in \mathbb{R}.
\end{equation}
    \end{enumerate}

    \item The following asymptotic development holds :
    \begin{equation}\label{int35}
\frac{E_0[1_{\Gamma_u}\varphi (S_t)]}{E_0[\varphi (S_t)]}=
Q_{0}^\varphi (\Gamma _u)+\sum_{i=1}^n\frac{1}{t^i}
E_0\big[1_{\Gamma_u}F_i^\varphi(X_u,S_u,u)\big] +
O\big(\frac{1}{t^{n+1}}\big),\quad t\rightarrow\infty.
\end{equation}
\end{enumerate}

\end{theo}

\noi Theorem \ref{tint7} will be proved  in Section \ref{asd}. We
also give a complement of Theorem \ref{tint7} (Theorem \ref{tasd1}
in Section \ref{asd}), taking as weight-process : $\psi
(S_t)e^{\lambda (S_t-X_t)}$, with $\lambda
>0$.

%

\section{Proof of Theorem \ref{tint2}, and of Theorem \ref{tint1}, as a consequence } \label{prt}
\setcounter{equation}{0}

 \noi Our proof of Theorem \ref{tint2} is
based on the following Lemma.

\begin{lemma} \label{lprt1} Let $y>0$, $u\geq 0$,  $\Gamma_u \in {\cal
F}_u$ and $t>u$. Then :
$$P_0(\Gamma_u | S _t=y)=\frac{p _{S _u}(y)}{p _{S _t}(y)}
E_0\Big[1_{\Gamma_u} h(t-u,y-X _u)|S _u=y\Big]\qquad \qquad$$
\begin{equation} \label{prt1}
\qquad \qquad \qquad \qquad \qquad \qquad \qquad+\frac{1}{p _{S
_t}(y)} E_0\Big[1_{\Gamma_u} 1 _{\{S_ u<y\}}p _{S
_{t-u}}(y-X_u)\Big].
\end{equation}

\noi where $p _{S _r}$ denotes the density function of $S_ r$, for
a fixed $ r>0$ :
\begin{equation} \label{prt2}
p _{S _r}(z)=\sqrt{ \frac{2} {\pi r} } e^{-z^2/2r} 1 _{\{z>0\}},
\end{equation}

\noi and
\begin{equation}\label{prt3}
h(r,z)=P(S _r<z)=\int _0^z p _{S _r}(x)dx=\sqrt{ \frac{2} {\pi r}
} \int _0^z e^{-x^2/2r}dx, \quad r,z>0.
\end{equation}
\end{lemma}

\begin{prooff} \ {\bf of Lemma \ref{lprt1}} Let $u\geq 0$,  $\Gamma_u \in {\cal
F}_u$ and $t>u$. It is clear that :
\begin{equation}\label{prt3a}
    S_t=S_u\vee\big(X_u +\max_{0\leq v\leq
t-u}\{X_{u+v}-X_u\}\big).
\end{equation}
Consequently if $g :[0,+\infty[ \rightarrow [0,+\infty]$ is Borel,
applying the Markov property at time $u$ leads to :
$$E_0\big[1_{\Gamma_u}g(S_t)\big]=E_0\big[1_{\Gamma_u}\tilde{g}(X_u,S_u)\big],$$
where
$$\tilde{g}(x,y)=E_0\big[g(y\vee\{x+S_{t-u}\})\big], \quad x_+ \leq y.$$
Then we easily obtain :
$$\begin{array}{ccl}
  \tilde{g}(x,y) & = & g(y)P_0(S_{t-u}\leq y-x)
  +E_0 \big[g(x+S_{t-u})1_{\{S_{t-u}>y-x\}}\big]
  \\
  &&\\
   & = &\dis g(y)h(t-u,y-x)+\int_0^\infty g(z)p_{S_{t-u}}(z-x)1_{\{z>y\}}dz. \\
\end{array}
$$
This proves (\ref{prt1}).
\end{prooff}


\begin{prooff} \ {\bf of Theorem \ref{tint2}} The two
estimates :
\begin{equation}\label{prt4}
    p_{S_t}(y) \sim \sqrt{\frac{2}{\pi t}}, \quad h(t,y)\sim
y\sqrt{\frac{2}{\pi t}}, \quad t\rightarrow \infty \ ( y>0),
\end{equation}
directly imply that $Q_{0,t}^{(y)}$  converges weakly to
$\widetilde{Q}_0^{(y)}$, as $t\rightarrow\infty$, where :
$$\widetilde{Q}^{(y)}_0(\Gamma_u)=p_{S_u}(y)E_0\Big[1_{\Gamma_u}(y-X_u)|S_u=y\Big] +
    E_0[1_{\Gamma_u}1_{\{S_u<y\}}], \quad \forall u\geq 0 \ \mbox{and}\ \Gamma_u\in {\cal F}_u.$$
Thanks to (\ref{int3}), $\dis 1-\Phi(y)=\int_y^\infty \varphi
(z)dz, y\geq 0$, then :
$$\begin{array}{ccl}
 \dis  \int _0 ^\infty  \widetilde{Q}^{(y)}_0(\Gamma_u)\varphi(y)dy & = & \dis E_0\Big[1_{\Gamma_u}\big(
(S_u-X_u)\varphi(S_u)+1-\Phi(S_u)\big)\Big] \\
   & = & E_0\big[ 1_{\Gamma_u}M_u^\varphi \big]=Q_0^\varphi(\Gamma_u).\\
\end{array}
$$
Consequently (\ref{int6}) implies $ \widetilde{Q}^{(y)}_0=
Q^{(y)}_0$, $ \ Q_0^\varphi (S_\infty \in dy)$  a.e.

\end{prooff}

\begin{rem}\label{rprt1}
     It is interesting to point out that (\ref{int9}) permits
    to prove that $y \mapsto Q_0^{(y)}$ is continuous, as the
    space of p.m.'s on the canonical space is endowed with the topology of weak
    convergence.
    \end{rem}

     \noi As indicated in Section \ref{int}, we now show how
     to prove Theorem \ref{tint1}, i.e. how to   recover
    (\ref{int4}) from Theorem \ref{tint2} and (\ref{int6}).

    \noi Indeed, let $\varphi$ be  as in Theorem \ref{tint1}. We have :
$$
 \frac{E_0\big[1_{\Gamma_u} \varphi (S_t)
\big] }{ E_0 \big[\varphi (S_t)\big]}=\frac{\int_0^\infty
 Q_{0,t}^{(y)}(\Gamma_u)\varphi(y) p_{S_t}(y)dy}{\int_0^\infty \varphi(y)
p_{S_t}(y)dy},
$$
where $u\geq 0, \Gamma_u \in{\cal F}_u$ and $t>u$.

\noi Using Theorem \ref{tint2}, (\ref{prt4}) and the dominated
convergence theorem, we get :
$$ \lim_{t \to \infty} \;\frac{E_0\big[1_{\Gamma_u} \varphi (S_t)
\big] }{ E_0 \big[\varphi (S_t)\big]}=\frac{\int_0^\infty
 Q_{0}^{(y)}(\Gamma_u)\varphi(y) dy}{\int_0^\infty \varphi(y)
 dy}=\int_0^\infty
 Q_{0}^{(y)}(\Gamma_u)\varphi(y) dy.$$


%
%

\section{Penalization for long Brownian bridges perturbed by their one-sided maximum }
\label{bri}

\setcounter{equation}{0}

 \noi
We keep the notation given in Sections \ref{int} and \ref{prt}.

\noi Let $Q_{0,t}^x$ be the law of the Brownian bridge started at
$0$, ending at $x$, with length $t$ :
\begin{equation}\label{bri1}
Q_{0,t}^x(\Gamma_t):=E_0[\Gamma_t |X_t=x], \quad \Gamma_t \in
{\cal F}_t.
\end{equation}
\noi (note the difference with the p.m. $Q_{0,t}^{(x)}$ defined in
(\ref{int7})).

 \noi Here, we make a simple remark concerning   the weak limit of
$Q_{0,t}^x$ as $t\rightarrow \infty$.

\noi Indeed, we observe that this limit is equal to the Wiener
measure $P_0$ : if $u\geq 0$ and $\Gamma_u \in {\cal F}_u$  then :
\begin{equation}\label{bri2}
    \lim_{t\rightarrow \infty}\ Q_{0,t}^x(\Gamma_u)=P_0(\Gamma_u),
\end{equation}
which follows from the fact that  $(X_s,0\leq s\leq u)$ under
$Q_{0,t}^x$, may be represented as $\dis
(B_s-\frac{s}{t}B_t+\frac{s}{t}x,,0\leq s\leq u)$, where $(B_s)$
is a Brownian motion started at $0$.

\noi The asymptotic study of long Brownian bridges penalized by
their one-sided maximum is more involved; in fact, we determine
the weak limit $Q_{0}^{a,y}$ of $Q_{0,t}^{a,y}$ as $t\rightarrow
\infty$, where $Q_{0,t}^{a,y}$ is the p.m. defined in
(\ref{int11}). The result is stated in Theorem \ref{tint3}.

\noi We proceed as for the proof of Theorem \ref{tint2}. We need
to generalize Lemma \ref{lprt1},  taking conditional expectations
with respect to $(S_t,X_t)$.

\begin{lemma}\label{lbri1} Let $a\in \mathbb{R}, y>a_+$, $u\geq 0$,
 $\Gamma_u \in {\cal F}_u$ and $t>u$. Then :
$$P_0(\Gamma_u | X_t=a,S _t=y)=\frac{p _{S _u}(y)}{p _{X_t,S_t}(a,y)}
E_0\Big[1_{\Gamma_u}\Big( \int_{ \mathbb{R}_+}
 p_{X_{t-u},S_{t-u}}(a-X_u,\xi)1_{\{\xi<y-X_u\}}d\xi\Big)\big| S_u=y\Big] $$
\begin{equation}\label{bri3}
\qquad \qquad \qquad +\frac{1}{p _{X_t,S_t}(a,y)}
E_0\Big[1_{\Gamma_u} 1 _{\{S_
u<y\}}p_{X_{t-u},S_{t-u}}(a-X_u,y-X_u)\Big].
\end{equation}
where $p _{X_v,S _v}$ denotes the density function of $(X_v,S_ v),
v>0$ :
\begin{equation}\label{bri4}
    p _{X_v,S _v}(a,y)=\sqrt{ \frac{2} {\pi v^3} }(2y-a) e^{-(2y-a)^2/2v} 1 _{\{y>a_+\}},
\end{equation}

\end{lemma}

\begin{proof} \ We imitate the proof of Lemma \ref{lprt1}.

\noi Let $g :[0,+\infty[\times [0,+\infty[ \rightarrow
[0,+\infty]$ be a Borel function. Thanks to (\ref{prt3a}), we have
:
$$E_0\big[1_{\Gamma_u}g(X_t,S_t)\big]=E_0\big[1_{\Gamma_u}\tilde{g}(X_u,S_u)\big],$$
where
$$\tilde{g}(a,y)=E_0\big[g(a+X_{t-u},y\vee\{a+S_{t-u}\})\big], \quad a_+ \leq y.$$
It follows :
$$\tilde{g}(a,y)=\tilde{g}_1(a,y)+\tilde{g}_2(a,y),$$
with :
$$ \tilde{g}_1(a,y)=E_0\big[g(a+X_{t-u},y)1_{\{S_{t-u}\leq
y-a\}}\big],$$
$$\tilde{g}_2(a,y)=E_0
\big[g(a+X_{t-u},a+S_{t-u})1_{\{S_{t-u}>y-a\}}\big].$$
\noi Since :
$$\tilde{g}_1(a,y)=\int_{\mathbb{R}\times \mathbb{R}_+} g(b,y)
p_{X_{t-u},S_{t-u}}(b-a,\xi)1_{\{\xi<y-a\}}db d\xi,$$
\noi and
$$\tilde{g}_2(a,y)=\int_{\mathbb{R}\times \mathbb{R}_+} g(b,z)
p_{X_{t-u},S_{t-u}}(b-a,z-a)1_{\{z>y\}}db dz,$$
then (\ref{bri3}) follows immediately.

\end{proof}

\begin{prooff} \ {\bf of Theorem \ref{tint3}}
Let  $u\geq 0$ and $\Gamma_u \in {\cal F}_u$.

\noi 1) Using (\ref{prt4}) and
\begin{equation}\label{bri5}
p_{X_{t},S_{t}}(a,y) \sim \sqrt{\frac{2}{\pi t^3}}(2y-a),\
t\rightarrow\infty, \quad y\geq a_+,
\end{equation}
we get :
$$\begin{array}{cl}
 \dis  \lim_{t\rightarrow
    \infty}Q_{0,t}^{a,y}(\Gamma_u)= & \dis
    \frac{p_{S_u}(y)}{2y-a}
    E_0\Big[1_{\Gamma_u}\Big(\int_{(a-X_u)_+}^{y-X_u}\big(2\xi-a+X_u\big)d\xi \Big)\big|S_u=y\Big]  \\
   &\dis  +\frac{1}{2y-a}E_0\Big[1_{\Gamma_u}1_{\{S_u<y\}}\big(2y-a-X_u\big)\Big]. \\
\end{array}
$$
The first integral in the right-hand side of the previous identity
may be computed, which yields :
$$Q_{0}^{a,y}(\Gamma_u)=p_{S_u}(y)\frac{y-a}{2y-a}E_0\big[1_{\Gamma_u}\big(y-X_u)|S_u=y\big]
+\frac{1}{2y-a}E_0\Big[1_{\Gamma_u}1_{\{S_u<y\}}\big(2y-a-X_u\big)\Big].$$
2) The relations (\ref{int9}) and (\ref{prt2}) imply :
$$\begin{array}{ccl}
  (2y-a)Q_{0}^{a,y}(\Gamma_u) & = & (y-a)\big\{Q_{0}^{(y)}(\Gamma_u)
-E_0\Big[1_{\Gamma_u}1_{\{S_u<y\}}\big]\big\}
+E_0\Big[1_{\Gamma_u}1_{\{S_u<y\}}\big(2y-a-X_u\big)\Big] \\
&&\\
   & = & (y-a)Q_{0}^{(y)}(\Gamma_u)
   +E_0\Big[1_{\Gamma_u}1_{\{S_u<y\}}\big(y-X_u\big)\Big]. \\
\end{array}
$$
Applying (\ref{int6}) with $\dis \varphi_y=\frac{1}{y}1_{[0,y]}$,
we get :
$$\int_0^y
Q_{0}^{(z)}(\Gamma_u)dz=yE_0\Big[1_{\Gamma_u}M_u^{\varphi_y}\big].$$
But $\dis \Phi_y(z):=\int_0^z\varphi_y(r)dr=\frac{z \wedge y}{y}$,
consequently :
$$ M_u^{\varphi_y} =(S_u-X_u)\varphi_y(S_u)+1- \Phi_y(S_u)=  \frac{y-X_u}{y}1_{\{S_u<y\}}.
$$
This proves (\ref{int13}).

\end{prooff}

\noi We now consider the Brownian bridge penalized by a function
of its one-sided maximum (cf Theorem \ref{tint4}).

\begin{prooff} \ {\bf of Theorem \ref{tint4}}

\noi Theorem \ref{tint4} is a direct consequence of Theorem
\ref{tint3}.

\noi Let  $u\geq 0$ and $\Gamma_u \in {\cal F}_u$ and $\varphi$ as
in Theorem \ref{tint4}.

\noi The relations (\ref{prt2}) and (\ref{bri4}) imply :
$$ P(S_t\in
dy|X_t=a)=\frac{2}{t}(2y-a)e^{-\frac{2y(y-a)}{t}}1_{\{y>a_+\}}
dy.$$
Consequently :
$$E_0\big[1_{\Gamma_u}\varphi(S_t)|X_t=a\big]=\frac{2}{t}
\int_0^{a_+}P_0(\Gamma_u|X_t=a,S_t=y)\varphi(y)(2y-a)e^{-\frac{2y(y-a)}{t}}
dy.$$
Hence :
$$\frac{E_0\big[1_{\Gamma_u}\varphi(S_t)|X_t=a\big]}{E_0\big[\varphi(S_t)|X_t=a\big]}
=\frac{\int_0^{a_+}P_0(\Gamma_u|X_t=a,S_t=y)\varphi(y)(2y-a)e^{-\frac{2y(y-a)}{t}}
dy}{\int_0^{a_+}\varphi(y)(2y-a)e^{-\frac{2y(y-a)}{t}} dy}.$$
Applying Theorem (\ref{tint3}) and the dominated convergence
theorem we get :
$$Q_0^{a,\varphi}(\Gamma_u)=\frac{\int_0^{a_+}(2y-a)\varphi(y)Q_0^{a,y}(\Gamma_u)
dy}{\int_0^{a_+}(2y-a)\varphi(y) dy}.$$
This proves (\ref{int16}). As for (\ref{int17}), it is a direct
consequence of (\ref{int13}).

%
%

%

\section{Proof of Theorem \ref{tint5} } \label{fXS}
\setcounter{equation}{0}

\noi 1) Point 1. of Theorem \ref{tint5} is a direct consequence of
Lemma \ref{lbri1} and Theorem \ref{tint3}.

\noi Taking the conditional expectation with respect to
$(X_t,S_t)$,  we obtain :
\begin{equation}\label{fXS1a}
\frac{E_0\big[1_{\Gamma_u} f(X_t, S_t) \big]}{E_0\big[ f(X_t, S_t)
\big]}=\frac{\int_\mathbb{R}da\int_{a_+}^\infty
Q_{0,t}^{a,y}(\Gamma_u)
p_{X_t,S_t}(a,y)f(a,y)dy}{\int_\mathbb{R}da\int_{a_+}^\infty
p_{X_t,S_t}(a,y)f(a,y)dy},
\end{equation}
where $p_{X_t,S_t}$ denotes the density function of $(X_t,S_t)$,
as given by (\ref{bri4}).

\noi Since $f$ satisfies (\ref{int18}), we may apply the dominated
convergence theorem; then taking the limit $t\rightarrow \infty$,
Theorem \ref{tint3}, and (\ref{bri5}) imply :
$$\lim_{t \to \infty} \; \frac{E_0\big[1_{\Gamma_u} f(X_t,
S_t) \big] }{ E_0
\big[f(X_t,S_t)\big]}:=\widetilde{Q}_0(\Gamma_u),$$
where :
$$\widetilde{Q}_0(\Gamma_u)=f^\star
\int_\mathbb{R}da\int_{a_+}^\infty
(2y-a)f(a,y)Q_0^{a,y}(\Gamma_u)dy.$$

\noi 2) We need to identify $\widetilde{Q}_0(\cdot)$.

\noi Let $\varphi$ be the function defined by (\ref{int20}) and
$\dis \Phi(y)=\int_0^y \varphi(z)dz, \ y \geq 0$.

\noi It is clear that $\varphi \geq 0$, then applying  Fubini's
theorem, we easily obtain :
\begin{equation}\label{fXS1}
\Phi(y)=f^\star\Big[\int_{\mathbb{R}\times\mathbb{R}_+}
f(a,\eta)1_{\{\eta >y\vee a_+\}}\big( \eta \wedge y
+(\eta-a)1_{\{\eta <y\}}\big)da d\eta \Big].
\end{equation}
In particular, taking the limit $y\rightarrow \infty$, we get :
$\dis \lim_{y\rightarrow \infty} \Phi(y)=1$. This means that
$\varphi$ satisfies (\ref{int3}).

\noi Moreover :
\begin{equation}\label{fXS2}
1-\Phi(y)=f^\star\Big[\int_{\mathbb{R}\times\mathbb{R}_+}
f(a,\eta)1_{\{\eta >y \vee a_+\}}\big( 2\eta -a-y)\big)da d\eta
\Big].
\end{equation}
Applying identity (\ref{int13}), we get :
$$\begin{array}{ccl}
  \widetilde{Q}_0(\cdot) & = & \dis f^\star
\int_\mathbb{R}da\Big(\int_{a_+}^\infty\Big\{(y-a)Q_0^{(y)}(\cdot)+\int_0^y
Q_0^{(\eta)}(\cdot)d \eta\Big\}f(a,y)dy\Big)\\
&  &  \\
   & = &\dis f^\star
\int_{\mathbb{R}\times \mathbb{R}_+}
\Big\{(\eta-a)f(a,\eta)1_{\{\eta>a_+\}}+\int_{\eta \vee
a_+}^\infty f(a,y)dy
\Big\}Q_0^{(\eta)}(\cdot) da d \eta \Big\}\\
&  &  \\
   & = & \dis \int_0^\infty \varphi(\eta)Q_0^{(\eta)}(\cdot) d\eta. \\
\end{array}
$$
Property (\ref{int6}) implies : $\widetilde{Q}_0=Q_0^{\varphi}$.

\noi 3) It remains to prove (\ref{int21}).

\noi Let $\widetilde{M}_t$ be the process defined as the
right-hand side of (\ref{int21}) :
$$
\widetilde{M}_t=f^\star \int_\mathbb{R}db\int_{
    b_+}^\infty (2y-b)f\big(b+X_t,S_t\vee(y+X_t)\big)dy.
$$
\noi Setting : $a=b+X_t$ and $\eta=y+X_t$, we obtain :
$$\widetilde{M}_t=f^\star
\int_\mathbb{R}\widetilde{\psi}_t(a)da,$$
where :
$$
\widetilde{\psi}_t(a)=\int_\mathbb{R}(2\eta-X_t-a)f(a,S_t\vee
\eta)1_{\{\eta >a , \eta >X_t\}}d\eta.$$
We have :
$$\begin{array}{ccl}
  \widetilde{\psi}_t(a) & = & \dis 1_{\{S_t>a\}}f(a,S_t)\int_{a\vee X_t}^{S_t}(2\eta-X_t-a)d\eta
  +\int_\mathbb{R}(2\eta-X_t-a)f(a, \eta)1_{\{\eta >S_t \vee a \}}d\eta
   \\
   &&\\
   & = &  \dis 1_{\{S_t>a\}}f(a,S_t)\big(S_t-a\vee
   X_t\big)\big((S_t+a\vee X_t -X_t-a\big)
  +\int_\mathbb{R}(2\eta-X_t-a)f(a, \eta)1_{\{\eta >S_t \vee a \}}d\eta\\
   & &  \\
   & =&\dis 1_{\{S_t>a\}}f(a,S_t)\big(S_t-
   X_t\big)\big((S_t-a\big)
  +\int_\mathbb{R}(2\eta-X_t-a)f(a, \eta)1_{\{\eta >S_t \vee a \}}d\eta.\\
\end{array}
$$
Since $M_t^\varphi =(S_t-X_t)\varphi (S_t)+1-\Phi(S_t)$, using
(\ref{int20}) and (\ref{fXS2}), we get :
$$M_t^\varphi=f^\star
\int_\mathbb{R}\psi_t(a)da,$$
where :
$$\psi_t(a)=(S_t-X_t)\Big[\int_\mathbb{R}f(a,\eta)1_{\{\eta >S_t\vee
    a_+\}} d\eta +f(a,S_t)(S_t-a)1_{\{a<S_t\}}\Big]
   +\int_\mathbb{R}f(a,\eta)(2\eta -a -S_t)1_{\{\eta >S_t\vee
    a_+\}}d\eta . $$
It is now clear that $\psi_t(a)=\widetilde{\psi}_t(a)$.
Consequently $M_t^\varphi=\widetilde{M}_t$.

\end{prooff}

\begin{rem}\label{rfXS1}\begin{enumerate}
    \item Let $f$ be of the type :
    $f(a,y)=f_1(a)1_{[0,A]}(y)$, where $A>0$ and $f_1 :]-\infty,A]
    \mapsto \mathbb{R}_+$ satisfies $\dis
    \int_{-\infty}^A(1+|a|)f_1(a)da<\infty$. Then it is easy to
    check that $\dis \bar f
    =A\int_{-\infty}^A(A-a)f_1(a)da<\infty$, and $\dis \varphi
    (y)=\frac{1}{A}1_{[0,A]}(y)$.

    \item It is possible to recover the identity (\ref{int13}) from Theorem
\ref{tint5}.

\noi Let $f$ as in Theorem \ref{tint5}. Using the first part of
the proof of Theorem \ref{tint5}, (\ref{int6}) and (\ref{int20}),
we have :
$$f^\star
\int_\mathbb{R}da\int_{a_+}^\infty
(2y-a)f(a,y)Q_0^{a,y}(\cdot)dy=Q_0^\varphi(\cdot)= \int_0^\infty
\varphi(y)Q_0^{(y)}(\cdot)dy,$$
$$\int_\mathbb{R}da\int_{a_+}^\infty
(2y-a)f(a,y)Q_0^{a,y}(\cdot)dy=\int_0^\infty
Q_0^{(y)}(\cdot)dy\Big\{\int_\mathbb{R}da \Big[\int_{y\vee
    a_+}^\infty f(a,\eta)d\eta +f(a,y)(y-a)1_{\{a<y\}}\Big]\Big\}.$$
Using Fubini's theorem, we easily obtain :
\begin{equation}\label{fXS3}
    \int_\mathbb{R}da\int_{a_+}^\infty
(2y-a)f(a,y)Q_0^{a,y}(\cdot)dy=\int_\mathbb{R}da\int_{a_+}^\infty
f(a,y)\Big[ y-a+\int_0^yQ_0^{(\eta)}(\cdot)d\eta\Big]dy,
\end{equation}
for any non-negative function $f$, satisfying (\ref{int18}), but
an easy application of Beppo-Levi theorem shows that (\ref{fXS3})
holds even without  (\ref{int18}) being satisfied.
\end{enumerate}
\end{rem}

%
%

\section{Penalization with $e^{ \lambda S_t+\mu X_t}$}
\label{pexp}

\setcounter{equation}{0}

1) In this section we focus on  penalizations with
weight-processes $f(X_t,S_t)$, where the function $f :
\mathbb{R}\times \mathbb{R}_+ \rightarrow \mathbb{R}_+$ belongs to
the family $\big\{f_{\lambda ,\mu}; \ f_{\lambda
,\mu}(a,y)=e^{\lambda y +\mu a}, \lambda , \mu \in
\mathbb{R}\big\}$.

\noi First, let us determine under which condition $f_{\lambda
,\mu}$ satisfies (\ref{int18}).

\noi Using the Fubini theorem, we have :
$$\overline{f}_{\lambda,\mu}=\int_0^\infty e^{\lambda y}dy
\int_{-\infty}^y (2y-a)e^{\mu a}da.$$
Consequently if $\mu \leq 0$ then
$\overline{f}_{\lambda,\mu}=\infty$.

\noi Suppose that $\mu >0$. The integral with respect to $da$ may
be computed, this yields to :
$$\overline{f}_{\lambda,\mu}=\frac{1}{\mu ^2}\int_0^\infty (1+\mu y)e^{(\lambda + \mu)y}dy.$$
As a result :
\begin{equation}\label{pexp1}
    \overline{f}_{\lambda,\mu}< \infty \quad \mbox{iff} \quad \mu
    >0 \ \mbox{ and } \ \lambda + \mu <0.
\end{equation}
Consequently if this condition holds, then Theorem \ref{tint5}
applies.

\noi 2) In our approach it is convenient to introduce $P_0^{\mu}$,
 the law of Brownian motion with drift $\mu$, starting at
$0$, and  $P_0^{(3)}$ the law of a three dimensional Bessel
process started at $0$.

\noi Recall    Pitman's theorem(\cite{P75}, \cite{J79}) :
\begin{enumerate}
    \item under $P_0$, the process $\big((2S_t-X_t,S_t), \ t\geq 0\big)$ is
distributed as $\big((X_t,J_t), \ t\geq 0\big)$ under $P_0^{(3)}$,
where $\dis J_t=\inf_{u\geq t}X_u$.
    \item let $\big( {\cal R}_t\big)$ be the natural filtration
    associated with the process $\big(R_t=2S_t-X_t, \ t\geq 0\big)$,
    then :
    \begin{equation}\label{pexp17}
    E_0[f(S_t)|{\cal
    R}_t]=\frac{1}{R_t}\int_0^{R_t}f(u)du,
\end{equation}
for any Borel function $f : \mathbb{R}_+\rightarrow \mathbb{R}_+$.

\end{enumerate}

\noi Pitman's theorem has been extended to the case of Brownian
motion with drift. From \cite{RP} (see also \cite{MY}), we know
that $(2S_t-X_t,t\geq 0)$ is a diffusion with generator :
\begin{equation}\label{pexp5}
    \frac{1}{2}\frac{d^2}{dx^2}+\mu \coth(\mu x)\frac{d}{dx}.
\end{equation}
%


\begin{prooff} \ {\bf of Theorem \ref{tint6}}


\noi Let $u$ be a fixed positive real number, $\Gamma_u \in {\cal
F}_u$, and define :
$$\Delta (\Gamma_u,t):=E_0\big[1_{\Gamma_u}e^{\lambda S_t +\mu
X_t}\big].$$

\noi 1) First suppose that $(\lambda , \mu)$ belongs to $ R_1$. We
have already proved that if $\mu >0$ then Theorem \ref{tint6} is a
direct consequence of Theorem \ref{tint5}. If $\mu=0$ and $\lambda
+\mu=\lambda <0$, then Theorem \ref{tint6} follows from Theorem
\ref{tint1}.

\noi 2) We now investigate the last case : $(\lambda , \mu)\in
R_3$.

\noi If $\mu=0$, then $\lambda <0$ and Theorem \ref{tint6} is a
direct consequence of Theorem \ref{tint1}.

\noi We suppose, in the sequel $\mu<0$.

 \noi We write $\Delta (\Gamma_u,t)$ as follows :
$$\Delta (\Gamma_u,t)=e^{\mu^2t/2}E_0^{\mu}
\big[1_{\Gamma_u}e^{\lambda S_t}\big].$$

\noi Applying the Markov property at time $u$, we get :
\begin{equation}\label{pexp6}
    \Delta (\Gamma_u,t)=e^{\mu^2t/2}E_0^{\mu}
\big[1_{\Gamma_u} h(X_u,S_u,t-u)\big],
\end{equation}
where
$$h(a,y,r)=E_0^{\mu}\big[e^{\lambda(y\vee (a+S_r))}\big],\quad
y\geq a_+, r\geq 0.$$
Since $\mu <0$, it is well-known that, under $P_0^\mu$,
$X_t\rightarrow-\infty$ as $t\rightarrow\infty$, $S_\infty<\infty$
and $P_0^\mu(S_\infty >x)=e^{2\mu x}, x\geq 0$.

\noi Consequently :
$$\lim_{r\rightarrow\infty}h(a,y,r)=I:=-2\mu \int_0^\infty e^{\lambda(y\vee (a+z))}e^{2\mu
z} dz.$$
Obviously the above integral  may be computed explicitly :
$$\begin{array}{ccl}
  I & = & \dis -2\mu \Big[e^{\lambda y}\int_0^{y-a} e^{2\mu z} dz+
 e^{\lambda a} \int_{y-a}^\infty e^{(\lambda +2\mu) z}dz \Big]\\
 &&\\
   & = & \dis  e^{\lambda y}\big[1-e^{2\mu(y-a)}+\frac{2\mu}{\lambda + 2 \mu}
e^{2\mu (y-a)}\Big]= e^{\lambda y}\big[1-\frac{\lambda}{\lambda +
2 \mu}e^{2\mu(y-a)}\Big].
   \\
\end{array}
$$
Moreover, it is easy to check :
$$e^{(\lambda +\mu)y-\mu a}\Big[\cosh\big(\mu(y-a)\big)
-\frac{\lambda +\mu}{\mu}\sinh\big(\mu(y-a)\big)\Big]=
\frac{\lambda +2\mu}{2\mu}e^{\lambda
y}\big[1-\frac{\lambda}{\lambda + 2 \mu}e^{2\mu(y-a)}\Big]$$
Finally :
$$\lim_{r\rightarrow\infty}h(a,y,r)=\frac{2\mu}{\lambda
+2\mu}e^{(\lambda +\mu)y-\mu a}\Big[\cosh\big(\mu(y-a)\big)
-\frac{\lambda +\mu}{\mu}\sinh\big(\mu(y-a)\big)\Big].$$
Coming back to (\ref{pexp6}), we obtain :
$$
\begin{array}{ccl}
\dis   \lim_{t \to \infty} \; \frac{E_0\big[1_{\Gamma_u} e^{\mu
X_t + \lambda S_t} \big] }{ E_0 \big[e^{\mu X_t + \lambda S_t}
\big]} & = &\dis \lim_{t \to \infty} \; \frac{E_0^{\mu}
\big[1_{\Gamma_u}h(X_u,S_u,t-u)\big]}{E_0^{\mu}
\big[h(0,0,t)\big]}
 \\
   & = & \dis E_0^{\mu}
\big[1_{\Gamma_u}e^{(\lambda + \mu)S_u-\mu X_u}\big\{
\cosh\big(\mu(S_u-X_u)\big) -\frac{\lambda
+\mu}{\mu}\sinh\big(\mu(S_u-X_u)\big)\big\}\Big]
\\
& = &\dis E_0 \big[1_{\Gamma_u}M^{\mu, \lambda}_u\big].
\end{array}
$$
%
 \noi 3) Let  $(\lambda , \mu)$ be an element of $ R_2$.

 \noi a) Let us start with the additional assumption : $\lambda + 2 \mu
 >0$. Since

%
$$P_0^{\lambda +\mu} =e^{(\lambda +\mu)X_t-(\lambda
+\mu)^2t/2} P_0  \quad \mbox{on } \ {\cal F}_t, $$
%
we have :
\begin{equation}\label{pexp2}
    \Delta (\Gamma_u,t)=e^{(\lambda+\mu)^2t/2}E_0^{\lambda +\mu}
\big[1_{\Gamma_u}e^{\lambda (S_t-X_t)}\big].
\end{equation}
Recall Theorem 1.1 in \cite{MY} : under $P_0^{\lambda +\mu}$, the
process $(S_t-X_t; t \geq 0)$ is distributed as $(|Y_t|, t\geq
0)$, where $(Y_t)$ is the so-called bang-bang process with
parameter $\lambda + \mu$, i.e. the diffusion with infinitesimal
generator :
\begin{equation}\label{pexp3}
    \frac{1}{2}\frac{d^2}{dx^2}-(\lambda +\mu) sgn(x)\frac{d}{dx}.
\end{equation}
\noi Applying the Markov property at time $u$ in (\ref{pexp2}),
yields to :
\begin{equation}\label{pexp4}
    \Delta (\Gamma_u,t)=e^{(\lambda+\mu)^2t/2}E_0^{\lambda +\mu}
\big[1_{\Gamma_u}\mathbb{E}_{S_u-X_u}\big\{e^{\lambda
|Y_{t-u}|}\big\}\big],
\end{equation}
where $\mathbb{P}_x$ denotes a p.m. under which $(Y_t)$ is the
diffusion process with generator (\ref{pexp3}) starting at $x$.

\noi Under  $\mathbb{P}_x$,  $(Y_t)$ is a recurrent diffusion  and
$\nu (dx):=(\lambda +\mu)e^{-2(\lambda +\mu)|x|}dx$ is its
invariant p.m.

\noi Consequently, for any $x \in \mathbb{R}$,
$$\lim_{r\rightarrow \infty}\mathbb{E}_x
\big[e^{\lambda |Y_r|}\big]=(\lambda +\mu)
\int_\mathbb{R}e^{\lambda |y|}e^{-2(\lambda
+\mu)|y|}dy=\frac{2(\lambda +\mu)}{\lambda +2\mu}.$$
Since  $\lambda + 2 \mu >0$ and $(\lambda , \mu)\in R_2$, then the
integral in the right-hand side is finite and does not depend on
$x$. As a result :
$$ \lim_{t \to \infty} \; \frac{E_0\big[1_{\Gamma_u} e^{\mu X_t +
\lambda S_t} \big] }{ E_0 \big[e^{\mu X_t + \lambda S_t} \big]}=
\lim_{t \to \infty} \; \frac{E_0^{\lambda +\mu}
\big[1_{\Gamma_u}\mathbb{E}_{S_u-X_u}\big\{e^{\lambda
|Y_{t-u}|}\big\}\big]}{\mathbb{E}_0 \big[e^{\lambda |Y_{t}|}\big]}
=P_0^{\lambda +\mu}(\Gamma_u).$$
Let us deal with the case $\lambda + 2 \mu=0, \ \mu \not=0$.
Applying (\ref{pexp2}), we have :
$$\Delta (\Gamma_u,t)=e^{\mu^2t/2}E_0^{-\mu}
\big[1_{\Gamma_u}e^{-\mu(2S_t-X_t)}\big].$$
The result follows from Pitman's theorem (for Brownian motion with
drift).

\noi b) It remains to study the case : $\lambda + 2 \mu =0$ and
$\lambda >0$.

\noi i) To begin with, we  modify $\Delta (\Gamma_u,t)$, $\lambda$
and $\mu$ being for now  two real numbers, without restriction.

\noi Applying the Markov property at time $u$ leads to :
\begin{equation}\label{pexp7}
    \Delta (\Gamma_u,t)=E_0
\big[1_{\Gamma_u} g(X_u,S_u,t-u)\big],
\end{equation}
where
$$g(a,y,r)=E_0\big[e^{\lambda\{y\vee (a+S_r)\}+\mu (a+X_r)}\big],\quad
y\geq a_+, r\geq 0.$$
Obviously, $g(a,y,r)$ may be decomposed  as follows :
\begin{equation}\label{pexp8}
g(a,y,r)=e^{\lambda y+\mu a}g_1(a,y,r)+e^{(\lambda+\mu)
a}g_2(a,y,r),
\end{equation}
with :
\begin{equation}\label{pexp9}
    g_1(a,y,r)=E_0\big[e^{\mu X_r} 1_{\{S_r<y-a\}}\big],\quad
    g_2(a,y,r)=E_0\big[e^{\lambda S_r+\mu X_r} 1_{\{S_r\geq
    y-a\}}\big].
\end{equation}

\noi Using Pitman's theorem recalled at the beginning of this
section, we get  :
\begin{equation}\label{pexp10}
     g_1(a,y,r)=  E_0^{(3)}\big[e^{-\mu X_r+2\mu J_r} 1_{\{J_r<y-a\}}\big]
  =E_0^{(3)}\big[e^{-\mu X_r}\frac{1}{X_r}\int_0^{(y-a)\wedge X_r}
   e^{2\mu z}dz\big],
\end{equation}
\begin{eqnarray}
     g_2(a,y,r)&= & E_0^{(3)}\big[e^{-\mu X_r+(\lambda +2\mu) J_r} 1_{\{J_r\geq y-a\}}\big]
     \nonumber\\
 & =&E_0^{(3)}\big[e^{-\mu X_r}1_{\{X_r\geq y-a\}}\frac{1}{X_r}\int_{y-a}^{ X_r}
   e^{(\lambda +2\mu) z}dz\big]. \label{pexp11}
\end{eqnarray}
As a result, if  $\mu \not= 0$  :
$$ g_1(a,y,r)\sim \frac{e^{2\mu(y-a)}-1}{2\mu}
E_0^{(3)}\big[e^{-\mu X_r}\frac{1}{X_r}\big],\quad
r\rightarrow\infty$$
Recall that :
\begin{equation}\label{pexp12}
P_0^{(3)}(X_r \in dz)=\sqrt{\frac{2}{\pi r^3}}z^2e^{-z^2/2r}
1_{\{z>0\}}dz.
\end{equation}
Then :
$$E_0^{(3)}\big[e^{-\mu
X_r}\frac{1}{X_r}\big]=\sqrt{\frac{2}{\pi r^3}}\int_0 ^\infty z
e^{-\mu z-z^2/2r}dz$$
Setting $\dis b =\frac{z + \mu r}{\sqrt{r}}$, we get :
$$E_0^{(3)}\big[e^{-\mu
X_r}\frac{1}{X_r}\big]=\sqrt{\frac{2}{\pi r }}e^{\mu
^2r/2}\int_{\mu \sqrt{r}} ^\infty (b-\mu \sqrt{r})e^{-b^2/2}d b.$$
It turns out that if $\mu <0$   :
\begin{equation}\label{pexp13}
    g_1(a,y,r)\sim \big(1-e^{2\mu(y-a)}\big) e^{\mu ^2r/2},\quad r\rightarrow
\infty.
\end{equation}
ii) We suppose now that $\lambda =-2\mu >0$.

\noi We need to determine the asymptotic behaviour of $g_2(a,y,r)$
as $r\rightarrow \infty$.

\noi Using  (\ref{pexp11}) and (\ref{pexp12}) we have :
$$\begin{array}{ccl}
  g_2(a,y,r) & = & \dis E_0^{(3)}\big[e^{-\mu X_r}1_{\{X_r\geq y-a\}}\frac{X_r-y+a}{X_r}
   \big] \\
   & = & \dis \sqrt{\frac{2}{\pi r^3}}\int_0 ^\infty z(z-y+a)
e^{-\mu z-z^2/2r}dz \\
&&\\
   &  =& \dis \sqrt{\frac{2}{\pi r }}e^{\mu
^2r/2}\int_{\mu \sqrt{r}} ^\infty (b-\mu \sqrt{r})(\sqrt{r}b-\mu r-y+a)e^{-b^2/2}d b. \\
\end{array}
   $$
As a result :
\begin{equation}\label{pexp14}
    g_2(a,y,r)\sim 2\mu ^2 r e^{\mu ^2r/2},\quad r\rightarrow
\infty.
\end{equation}
Due to (\ref{pexp7}), (\ref{pexp8}), (\ref{pexp13}),
(\ref{pexp14}) and $\lambda=-2\mu$, we get :
$$
\Delta (\Gamma_u,t)\sim 2\mu ^2 t e^{\mu ^2t/2} E_0
\big[1_{\Gamma_u}e^{-\mu X_u-\mu^2 u/2}\big], \quad t\rightarrow
\infty.
$$
In particular :
$$\Delta (\Omega,t)=E_0 \big[e^{-\mu X_{t}+\lambda S_t}\big] \sim 2\mu ^2 t e^{\mu ^2t/2}
, \quad t\rightarrow \infty.
$$
Finally :
$$ \lim_{t \to \infty} \; \frac{E_0\big[1_{\Gamma_u} e^{\mu X_t +
\lambda S_t} \big] }{ E_0 \big[e^{\mu X_t + \lambda S_t} \big]}=
E_0 \big[1_{\Gamma_u}e^{-\mu X_u-\mu^2 u/2}\big].$$

\end{prooff}

\begin{rem} \label{rexp1} Here is another proof of Theorem
\ref{tint6} : keeping the notations introduced in point 3) b ) i)
of the proof above, recall that we have proved :
$$E_0\big[1_{\Gamma_u} e^{\mu X_t +
\lambda S_t} \big]=E_0\big[1_{\Gamma_u}\big\{e^{\mu X_u + \lambda
S_u}g_1(X_u,S_u,t-u)+e^{(\lambda +\mu) X_u
}g_2(X_u,S_u,t-u)\big\}\big],$$
where the functions $g_1(a,y,r)$ and $g_2(a,y,r)$ are given by
(\ref{pexp9}) or (\ref{pexp10}) and (\ref{pexp11}). In our proof
of Theorem \ref{tint6} we  only need the asymptotics of
$g_i(a,y,r)$ as $r\rightarrow \infty,\  i=1,2$ in the case
$\lambda + 2 \mu =0, \ \lambda >0$. It is actually possible to
determine the asymptotics of the previous quantities in any case.
However tedious calculations are needed, this explains why we have
given a short and direct proof of Theorem \ref{tint6}.

\end{rem}

 \noi We now give a direct interpretation of Theorem \ref{tint6}
in terms of  the three dimensional Bessel process and its
post-minimum.

\begin{prop} \label{ppexp1} Let $\lambda , \mu \in \mathbb{R}$.
\begin{enumerate}
    \item For every  $u \ge 0$, and $\Gamma_u$ in ${\cal F}_u$,
    \begin{equation}\label{pexp15}
   \lim_{t \to \infty} \; \frac{E_0^{(3)}\big[1_{\Gamma_u} e^{\mu X_t +
\lambda J_t} \big] }{ E_0^{(3)} \big[e^{\mu X_t + \lambda J_t}
\big]}:=E_0^{(3)}\big[1_{\Gamma_u}\overline{M}_u^{\mu, \lambda
}\big] ,
\end{equation}
 where   $(\overline{M}_u^{\mu, \lambda})$ is the positive $\big(({\cal F}_u),P_0^{(3)}\big)$
    martingale :
\begin{equation}\label{pexp16}
\overline{M}_u^{\mu, \lambda}= \left\{
\begin{array}{ll}
1 &  \mbox{if}
  \ \lambda + \mu  <0 \ \mbox{ and } \mu \leq 0,\\
  &\\
\dis
  e^{\{ -(\lambda + \mu)^2u/2\}}\frac{ \sinh \big((\lambda +\mu)X_u\big)}{(\lambda +\mu)X_u}
 \big]& \mbox{if}
  \ \lambda \geq 0 \ \mbox{ and } \lambda +\mu \geq 0, \\
  &\\
  \dis
   e^{ - \mu^2u/2}\frac{ \sinh \big(\mu X_u \big)}{\mu X_u}
 & \mbox{if}
  \ \lambda <0 \ \mbox{ and } \mu >0, \\
\end{array}
\right.
\end{equation}
Note that $\overline{M}_0^{\mu, \lambda}=1$.

\item The map :  $\Gamma_u (\in {\cal F}_u) \mapsto
E_0^{(3)}\big[1_{\Gamma_u}\overline{M}_u^{\mu, \lambda }\big]$
induces a p.m. on $\big(\Omega,{\cal F}_\infty\big)$.
\end{enumerate}
\end{prop}

\begin{proof} \ Proposition \ref{ppexp1} is a direct consequence of
Theorem \ref{tint6} and Pitman's theorem.

\noi We have :
$$\frac{E_0^{(3)}\big[1_{\Gamma_u} e^{\mu X_t +
\lambda J_t} \big] }{ E_0^{(3)} \big[e^{\mu X_t + \lambda J_t}
\big]}=\frac{E_0\big[1_{\widehat{\Gamma}_u} e^{-\mu X_t + (\lambda
+2\mu )S_t} \big] }{ E_0 \big[e^{-\mu X_t + (\lambda +2\mu )S_t}
\big]}$$
where $\widehat{\Gamma}_u:=\{\omega \in \Omega; \widehat{\omega}
\in \Gamma_u\}$, and $\dis \widehat{\omega}_t:=\sup_{0\leq u \leq
t}\omega (u)-\omega (t)$.

\noi Applying our Theorem \ref{tint6}, we obtain :
$$\lim_{t \to \infty} \;\frac{E_0^{(3)}\big[1_{\Gamma_u} e^{\mu X_t +
\lambda J_t} \big] }{ E_0^{(3)} \big[e^{\mu X_t + \lambda J_t}
\big]}=E_0\big[1_{\widehat{\Gamma}_u}M_u^{-\mu, \lambda +2\mu
}\big].$$
Since $\widehat{\Gamma}_u \in {\cal R}_u$ then :
$$E_0\big[1_{\widehat{\Gamma}_u}M_u^{-\mu, \lambda +2\mu
}\big]=E_0\Big[1_{\widehat{\Gamma}_u}E_0\big[M_u^{-\mu, \lambda
+2\mu }|{\cal R}_u\big]\Big].$$
We claim that :
\begin{equation}\label{pexp18}
E_0\big[M_u^{-\mu, \lambda +2\mu }|{\cal R}_u\big]= \left\{
\begin{array}{ll}
1 &  \mbox{if}
  \ \lambda + \mu  <0 \ \mbox{ and } \mu \leq 0,\\
  &\\
\dis
  e^{\{ -(\lambda + \mu)^2u/2\}}\frac{ \sinh \big((\lambda +\mu)(2S_u-X_u)\big)}
  {(\lambda +\mu)(2S_u-X_u)}
 & \mbox{if}
  \ \lambda \geq 0 \ \mbox{ and } \lambda +\mu \geq 0, \\
  &\\
  \dis
   e^{ - \mu^2u/2}\frac{ \sinh \big(\mu (2S_u-X_u) \big)}{\mu (2S_u-X_u)}
 & \mbox{if}
  \ \lambda <0 \ \mbox{ and } \mu >0, \\
\end{array}
\right.
\end{equation}
Making again use of Pitman's theorem, it  is immediate to obtain
(\ref{pexp16}).

\noi As for (\ref{pexp18}), we only prove the third case. The two
other cases may be proved similarly. Note that $(\lambda +2 \mu,
-\mu)\in R_1$ (resp. $R_2$) iff $\lambda + \mu < 0$ and $ \mu \leq
0$ (resp. $\lambda \geq 0$ and $\lambda + \mu \geq 0$).

\noi As for the third case, we have : $(\lambda +2 \mu, -\mu)\in
R_3$ iff $\lambda< 0$ and $ \mu < 0$.

\noi Setting $R_u:=2S_u-X_u$, then (\ref{int28})and (\ref{pexp17})
imply :
$$M_u^{-\mu, \lambda +2\mu }=e^{\{ (\lambda + \mu)S_u-
\mu^2u/2\}}\big[\cosh\big(\mu(R_u-S_u)\big)
  - \frac{\lambda + \mu}{\mu}\sinh\big(\mu(R_u-S_u)\big)\big],$$
$$\begin{array}{ccl}
  E_0\big[M_u^{-\mu, \lambda +2\mu }|{\cal R}_u\big] & = &
 \dis \frac{e^{ -\mu^2u/2}}{R_u}\int_0^{R_u}e^{ (\lambda + \mu)y}
 \big[\cosh\big(\mu(R_u-y)\big)
  - \frac{\lambda + \mu}{\mu}\sinh\big(\mu(R_u-y)\big)\big]dy \\
  &&\\
   & = & \dis \frac{e^{- \mu^2u/2}}{R_u}
   \Big[-\frac{1}{ \mu} e^{ (\lambda + \mu)y}
   \sinh\big(\mu(R_u-y)\big)\Big]_{y=0}^{y=R_u}\\
   &&\\
   & = & \dis e^{ -\mu^2u/2}\frac{\sinh\big(\mu R_u \big)}{\mu R_u}. \\
\end{array}
$$
This establishes the third case in (\ref{pexp16}), using again
Pitman's theorem.
\end{proof}

\begin{rem} \label{rexp2} It seems natural to ask for :
\begin{equation}\label{pexp18a}
    \lim_{t\rightarrow \infty}\frac{E_0^{(3)}\big[1_{\Gamma_u} f( X_t
, J_t) \big] }{ E_0^{(3)} \big[f(X_t , J_t) \big]},
\end{equation}
for some suitable Borel $f : \mathbb{R}_+ \times \mathbb{R}_+
\mapsto \mathbb{R}_+$.

\noi Using Pitman's theorem (see 1. in the proof of Proposition
\ref{ppexp1}), the above ratio  is equal to :
$$\frac{E_0\big[1_{\widehat{\Gamma}_u}f( 2S_t-X_t ,S_t) \big] }
{E_0 \big[f( 2S_t-X_t ,S_t) \big]}.$$
Consequently Theorem \ref{tint5} applies as soon as :
\begin{equation}\label{pexp19}
    \widetilde{f}:=\int_\mathbb{R} da \int_{a_+}^\infty (2y-a) f(2y-a,y)dy
=\int_{\mathbb{R}_+\times \mathbb{R}_+} f(b,y)1_{\{b>y\}}dbdy <
\infty.
\end{equation}
Suppose that this condition holds. Then
$$\lim_{t\rightarrow \infty}\frac{E_0^{(3)}\big[1_{\Gamma_u} f( X_t
, J_t) \big] }{ E_0^{(3)} \big[f(X_t , J_t)
\big]}=E_0\big[1_{\widehat{\Gamma}_u}M_u^{\varphi}\big],$$
with
$$\varphi(y)=f^\dagger \Big[
\int_{\mathbb{R}_+\times \mathbb{R}_+} f(b,\eta)1_{\{b>\eta
>y\}}dbd\eta +\int_y^\infty f(b,y)db \Big],$$
and $f^\dagger=1/\widetilde{f}$.

\noi Proceeding as in the proof of Proposition \ref{ppexp1} we may
prove :
\begin{equation}\label{pexp20}
E_0\big[M_u^{\varphi}| {\cal R}_u\big]=1.
\end{equation}
Finally the limit in (\ref{pexp18a}) equals $P_0^{(3)}(\Gamma
_u)$. In other words the penalization with $f(X_t,J_t)$, $f$
satisfying
 (\ref{pexp19}) does not generate a new p.m.

\end{rem}

\noi As an end to this section, we would like to discuss the
relationship between Theorem \ref{tint6} and the results obtained
in  \cite{HY}. Recall that these authors have proved that
\begin{equation}\label{pexp21}
\lim_{t\rightarrow \infty}\frac{E_0^{\mu}\big[1_{\Gamma_u}
A_t^{\lambda /2}  \big] }{ E_0^{\mu} \big[A_t^{\lambda /2} \big]},
\end{equation}
exists where $\lambda, \mu \in \mathbb{R}$, $P^\mu_0$ denotes the
p.m. on canonical space which makes $(X_t)$ a Brownian motion with
drift $\mu$, started at $0$, and :
$$A_t=\int_0^t e^{2X_s}ds, \quad t \geq 0.$$
As for our Theorem \ref{tint6},  it is proved in \cite{HY}, that
 a phase transition phenomenon occurs : there exists three
disjoint regions in $\mathbb{R}\times \mathbb{R}$ associated with
three types of limit distributions in (\ref{pexp21}). It is
striking to note that these regions coincide with the domains
$R_1,R_2$ and $R_3$ introduced in (\ref{int24})-(\ref{int26}).

\noi We have actually no proof of this fact. Nevertheless if
$\lambda >0$, we have a heuristic argument :

$$\frac{E_0^{\mu}\big[1_{\Gamma_u}
A_t^{\lambda /2}  \big] }{ E_0^{\mu} \big[A_t^{\lambda /2} \big]}=
\frac{E_0\big[1_{\Gamma_u} e^{\mu X_t} A_t^{\lambda /2}  \big] }{
E_0 \big[e^{\mu X_t} A_t^{\lambda /2} \big]}.$$
\noi Roughly speaking, the Laplace theorem tells us that  $\dis
A_t=\int_0^t e^{2X_s}ds$ has the same behaviour as $e^{2S_t}$, see
more precisely, the limit results in \cite{CMY} (formulae (61) and
(62) p 181) and \cite{PY96}. Therefore replacing formally $A_t$ by
$e^{2S_t}$, we get :
$$\frac{E_0\big[1_{\Gamma_u} e^{\mu X_t} A_t^{\lambda /2}  \big] }{
E_0 \big[e^{\mu X_t} A_t^{\lambda /2} \big]}
 \approx \frac{E_0\big[1_{\Gamma_u} e^{\mu X_t + \lambda S_t} \big] }{ E_0
\big[e^{\mu X_t + \lambda S_t} \big]},\quad  t\rightarrow
\infty.$$
%

%
\section{ Asymptotic development}\label{asd}
%

\setcounter{equation}{0} \noi We first recall  a penalization
result obtained in (\cite{RVY2a}), choosing as weight-process :
$\psi (S_t)e^{\lambda (S_t-X_t)}$, where $\lambda >0$.

\noi Let us start with some notations. Let $\psi : \mathbb{R}_+
\mapsto \mathbb{R}_+$ be a Borel function satisfying :
\begin{equation}\label{asd1b}
    \int_0^\infty \psi(z)e^{-\lambda z}dz=1.$$
\end{equation}
To $\psi$ we associate   the functions $\Phi$ and $\varphi$ :
\begin{equation}\label{asd1c}
    \Phi(y):=1-e^{\lambda y}\int _y^\infty \psi(z)e^{-\lambda
    z}dz, \quad y \geq 0,
\end{equation}
\begin{equation}\label{asd1d}
    \varphi (y):=\Phi'(y)=\psi (y)-\lambda e^{\lambda y}\int _y^\infty \psi(z)e^{-\lambda
    z}dz \quad y \geq 0.
\end{equation}
Then :
\begin{equation}\label{asd1e}
    M_t^{\lambda ,\varphi}:=\Big\{\psi (S_t)\frac{\sinh \big(\lambda (S_t-X_t)\big)}{\lambda}
    +e^{\lambda X_t}\int_{S_t}^\infty \psi(z)e^{-\lambda
    z}dz\Big\}e^{-\lambda ^2 t/2} \quad t \geq 0,
\end{equation}
is a $\big(({\cal F}_t),P_0\big)$ positive, and continuous
 martingale.

\noi In this setting we have proved (see Theorem 3.9 in
\cite{RVY2a}).

\begin{prop}\label{pasd1} Let $\psi : \mathbb{R}_+ \mapsto
\mathbb{R}_+$ be a Borel function satisfying (\ref{asd1b}). Then :
\begin{equation}\label{asd1f}
    \lim_{t\rightarrow \infty}\frac{E_0\Big[1_{\Gamma_u}\psi (S_t)
    e^{\lambda (S_t -X_t)}\Big]}{E_0\Big[\psi (S_t)
    e^{\lambda (S_t -X_t)}\Big]}=E_0[1_{\Gamma_u}M_u^{\lambda
    ,\varphi}\Big],
\end{equation}
for any $u\geq 0$ and $\Gamma_u \in {\cal F}_u$.

\end{prop}

\begin{rem} \label{rasd1}\begin{enumerate}
    \item In \cite{RVY2a}, we have determined the law of $(X_t)$ under
    the new p.m. $\Gamma_u (\in {\cal F}_u) \mapsto E_0\big[1_{\Gamma_u}M_u^{\lambda
    ,\varphi}\big]$. However, this result is not used in the
    sequel.
    \item If we take $\lambda =0$ and $\psi=\varphi$, then
    (\ref{asd1d}) holds, (\ref{asd1b}) corresponds to (\ref{int3})
    and
    $(M_t^{0 ,\varphi})$ coincides with the martingale
    $(M_t^{\varphi})$ defined by (\ref{int5}). With these
    conventions, Proposition \ref{pasd1} is an extension of points 1. and 2. of Theorem
    \ref{tint1}.
\end{enumerate}

\end{rem}

\noi The aim of this section is to prove that, under suitable
assumptions, we can obtain an asymptotic expansion of $\dis
t\mapsto \frac{E_0\Big[1_{\Gamma_u}\psi (S_t)
    e^{\lambda (S_t -X_t)}\Big]}{E_0\Big[\psi (S_t)
    e^{\lambda (S_t -X_t)}\Big]}$ as $t\rightarrow\infty$.
Note that Proposition \ref{pasd1} gives the first term.

\begin{theo}\label{tasd1} Let $\psi : \mathbb{R}_+ \mapsto
\mathbb{R}_+$  satisfying (\ref{asd1b}). We suppose that there
exists an integer $n\geq 1$ such that :
\begin{equation}\label{asd1g}
    \int_0^\infty \psi(y)(1+y^n)dy<\infty .
\end{equation}
\begin{enumerate}
    \item There exists a family of functions $(F_i^{\lambda,\varphi})_{1\leq i \leq
    n}$, $F_i^{\lambda,\varphi}  :  \mathbb{R}\times  \mathbb{R}_+\times
    \mathbb{R}_+ \mapsto \mathbb{R}$, such that
    \begin{enumerate}
        \item $(F_i^{\lambda,\varphi}(X_t,S_t,t), t\geq 0)$ is a $\big(({\cal
        F}_t),P_0 \big)$-martingale, for any $1\leq i \leq n$,
        \item If $i=1$, we have :
        \begin{equation}\label{asd1h}
F_1^{\lambda,\varphi}(X_t,S_t,t)= \frac{c(\lambda
,\varphi)}{\lambda^3\sqrt{2\pi}}\Big(M^{\varphi_1}_t
-M_t^{\lambda,\varphi}\Big).
$$
\end{equation}
%
where $\dis c(\lambda, \varphi):=\int_{0}^\infty \psi(x)(1
-\lambda x)dx$ and
\begin{equation}\label{asd1i}
   \varphi_1(y):=\frac{1}{c(\lambda, \varphi)}
   \Big(\psi(y)-\lambda \int_y^\infty \psi (x)dx\Big), \quad y\geq 0.
\end{equation}
(Note that  $\dis \int_0^\infty \varphi_1 (y)dy=1$).
    \end{enumerate}

    \item The following asymptotic development as $t\rightarrow\infty$, holds :
    \begin{equation}\label{asd1j}
\frac{E_0\big[1_{\Gamma_u}\psi (S_t)
    e^{\lambda (S_t -X_t)}\big]}{E_0\big[\psi (S_t)
    e^{\lambda (S_t -X_t)}\big]}= E_0\big[1_{\Gamma_u}M_u^{\lambda
    ,\varphi}\big]+\frac{e^{-\lambda ^2 t/2}}{\sqrt{t}}\Big(\sum_{i=1}^n\frac{1}{t^i}
E_0\big[1_{\Gamma_u}F_i^{\lambda,\varphi}(X_u,S_u,u)\big] +
O\big(\frac{1}{t^{n+1}}\big)\Big).
\end{equation}
\end{enumerate}

\end{theo}

\noi  Note that the two asymptotic expansions (\ref{int35}) and
(\ref{asd1j}) are drastically different, depending on whether
$\lambda
>0$ or $\lambda =0$. We have already observed in Remark
\ref{rasd1}, that taking formally $\lambda=0$ in (\ref{asd1f})
gives (\ref{int4}). In other words the first term in (\ref{asd1j})
(with $\lambda =0$) coincides with the first term in
(\ref{int35}). However the expansion is expressed in terms of
powers $t^{-(i+1/2)}$ instead of $t^{-i}$.

\begin{prooff} \ {\bf of Theorems \ref{tint7} and \ref{tasd1}}

\noi 1) Let us start with some common features concerning the two
cases $\lambda >0$ and $\lambda =0$, i.e. $\lambda \geq 0$. We
adopt the convention that $\psi=\varphi$ if $\lambda =0$.

\noi Let  $u$ be a fixed positive real number, $\Gamma_u \in {\cal
F}_u$, and
$$\Delta (\lambda,\Gamma_u,t):=E_0\big[1_{\Gamma_u}\psi (S_t)e^{\lambda ( S_t
-X_t)}\big],$$
where $\lambda \geq 0$.

 \noi

\noi Applying the Markov property at time $u$, we get :
\begin{equation}\label{asd1}
    \Delta (\lambda ,\Gamma_u,t)=E_0
\big[1_{\Gamma_u} g(\lambda ,X_u,S_u,t-u)\big],
\end{equation}
where
$$g(\lambda, a,y,r):=E_0\big[\psi\big(y\vee (a+S_r)\big)
e^{\lambda(y\vee (a+S_r)-a-X_r)}\big], \quad r \geq 0, y\geq a_+,
a\in \mathbb{R}.$$
We have :
$$g(\lambda, a,y,r)=\psi (y)e^{\lambda
(y-a)}E_0\big[e^{-\lambda X_r}1_{\{S_r<y-a\}}\big]+
E_0\big[\psi(a+S_r)e^{\lambda (S_r- X_r)}1_{\{S_r\geq
y-a\}}\big].$$
Using Pitman's theorem, we get :
%
\begin{eqnarray}
 g(\lambda, a,y,r) & = & \psi (y)e^{\lambda
(y-a)}E_0^{(3)}\big[e^{-\lambda(2J_r- X_r)}1_{\{J_r<y-a\}}\big]+
E_0^{(3)}\big[\psi(a+J_r)e^{\lambda (X_r- J_r)}1_{\{J_r\geq
y-a\}}\big] \nonumber\\
&&\nonumber\\
   & = & \dis E_0^{(3)}\Big[\frac{e^{\lambda X_r}}{X_r}\int_0^{X_r}
\big\{\psi (y)e^{\lambda (y-a)}e^{-2\lambda z}1_{\{z<y-a\}}+\psi
(a+z)e^{-\lambda z} 1_{\{z\geq
y-a\}}\big\}dz\Big]\nonumber\\
&&\nonumber\\
   & = & \dis \psi (y)e^{\lambda (y-a)}\int_0^{y-a}e^{-2\lambda
z}h(\lambda,z,r)dz+\int_{y-a}^\infty e^{-\lambda
z}\psi(a+z)h(\lambda,z,r)dz, \label{asd1a}
\end{eqnarray}
where :
$$h(\lambda,z,r)=E_0^{(3)}\big[\frac{e^{\lambda
X_r}}{X_r}1_{\{X_r>z\}}\big].
$$
%
%
%
Applying (\ref{pexp12}), we get  :
\begin{equation}\label{asd3}
h(\lambda,z,r)=\sqrt{\frac{2}{\pi r^3}}e^{\lambda ^2 r/2}\int
_z^\infty b e^{-(b-\lambda r)^2/2r}db=\sqrt{\frac{2}{\pi r}}
e^{\lambda ^2 r/2}\int _{\frac{z-\lambda r}{\sqrt{r}}}^\infty
(\lambda \sqrt{r}+v) e^{-v^2/2}dv.
\end{equation}

\noi 2) Suppose that $\lambda =0$. Therefore we replace in the
sequel  $\psi$ by $\varphi$.

\noi a) Then :
$$h(0,z,r)=\sqrt{\frac{2}{\pi r^3}}\int _z^\infty b e^{-b^2/2r}db
=\sqrt{\frac{2}{\pi r}}e^{-z^2/2r},$$
%

%
$$g(0, a,y,r)=\sqrt{\frac{2}{\pi r}}\widehat{g}(0, a,y,r),$$
with :
$$ \ \widehat{g}(0, a,y,r)= \varphi
(y)\int_0^{y-a}e^{-z^2/2r}dz+\int_{y}^\infty
e^{-(v-a)^2/2r}\varphi(v)dv.$$
Let us introduce :
$$A_i(a,y):=\varphi (y)\frac{(y-a)^{2i+1}}{2i+1}
+\int_{y}^\infty (v-a)^{2i}\varphi(v)dv,$$ (note that $\dis
A_i(0,0)= \int_0^\infty v^{2i}\varphi(v)dv$; in particular
$A_0(0,0)=1)$.

\noi Using the  series development of $e^{-\theta}$ with $\theta
\geq 0$, we get :
\begin{equation}\label{asd5}
     \widehat{g}(0, a,y,r)  = \sum_{i=0}^n \frac{(-1)^i}{(2r)^i
i!}A_i(a,y) +O\big(\frac{1}{r^{n+1}} \big).
\end{equation}
%
%

\noi Moreover $\varepsilon \mapsto
\widehat{g}(0,a,y,1/\varepsilon)$ is of class $C^\infty$ on
$[0,1/2]$ and :
\begin{equation}\label{asd6}
    |\frac{\partial^i\widehat{g}(0,a,y,1/\varepsilon)}{\partial \varepsilon^i}|
    \leq k_i \Big(\varphi (y)(y-a)^{2i+1}
+\int_{0}^\infty (v-a)^{2i}\varphi(v)dv\Big).
\end{equation}

\noi Suppose that $t \rightarrow \infty$, then :
\begin{eqnarray}
  \frac{\dis g(0, a,y,t-u)}{\dis g(0,0,0,t)} & = & \dis \big(1-\frac{u}{t}\big)^{-1/2}
 \  \frac{\dis \sum_{i=0}^n
\frac{(-1)^i}{2^i
i!}\frac{1}{t^i(1-u/t)^i}A_i(a,y)+O\big(\frac{1}{t^{n+1}}\big)}
{\dis \sum_{i=0}^n \frac{(-1)^i}{2^i
i!}\frac{1}{t^i}A_i(0,0)+O\big(\frac{1}{t^{n+1}}\big)} \nonumber \\
 & = & \dis \sum_{i=0}^n
\frac{1}{t^i}F_i(a,y,u)+\frac{1}{t^{n+1}}R(0,a,y,u,t),
\label{asd7}
\end{eqnarray}
where $t\mapsto R(0,a,y,u,t)$ is bounded, and $F_i(a,y,u)$ may be
written in the following form :
\begin{equation}\label{asd7b}
F_i(a,y,u)=\sum_{j=0}^i \alpha_{i,j}(u)A_j(a,y),
\end{equation}
$\alpha_{i,j}(u) $ being some polynomial function.

 \noi  b) In
particular :
\begin{equation}\label{asd8}
    F_0(a,y,u)=\frac{A_0(a,y)}{A_0(0,0)}=A_0(a,y)=\varphi (y)(y-a)
+\int_{y}^\infty \varphi(v)dv=\varphi (y)(y-a)+\Phi (y),
\end{equation}
(note that $F_0$ does not depend on $u$).

\noi To compute $F_1(a,y,u)$ we need the first order term   :
$$\frac{g(0,
a,y,t-u)}{g(0,0,0,t)}=\frac{\dis
A_0(a,y)\Big(1-\frac{1}{2t}\frac{A_1(a,y)}{A_0(a,y)}\Big)\Big(1-\frac{u}{t}\Big)^{-1/2}}
{\dis
A_0(0,0)\Big(1-\frac{1}{2t}\frac{A_1(0,0)}{A_0(0,0)}\Big)}+O\big(\frac{1}{t^2}\big).$$
Recall that $A_0(0,0)=1$, consequently :
$$\begin{array}{ccl}
  F_1(a,y,u) & = & \dis \frac{A_0(a,y)}{2}\Big(A_1(0,0)
-\frac{A_1(a,y)}{A_0(a,y)}+u\Big)\\
&&\\
   & = & \dis \frac{A_1(0,0)+u}{2}A_0(a,y)-\frac{A_1(a,y)}{2} \\
   &&\\
   & = & \dis  \frac{\int_0^\infty v^2\varphi (v)dv +u}{2}A_0(a,y)
   -\varphi (y)\frac{(y-a)^{3}}{3!}-\frac{1}{2}
\int_{y}^\infty (v-a)^{2}\varphi(v)dv.\\
\end{array}
$$

\noi c) We would like to obtain some estimates about the remainder
term $R(0,a,y,u,t)$ in (\ref{asd7}), as a function of $(a,y)$.

\noi Taking $t \geq 2u+2$ and setting $\varepsilon=1/t$ we have :
$\varepsilon \leq 1/2\ $, $\dis \varepsilon u\leq 1/2\ ,
\frac{1}{t-u}=\frac{\varepsilon}{1-u\varepsilon}\in [0,1]$

\noi Let
$\widehat{g}_1(0,a,y,\varepsilon):=\widehat{g}(0,a,y,(1-u\varepsilon)/\varepsilon),
\ \varepsilon \in ]0,1/(2u+2)]$. Then property  (\ref{asd6})
implies :
\begin{equation}\label{asd9}
|\frac{\partial^{i}\widehat{g}_1(0,a,y,\varepsilon)}{\partial
\varepsilon^{i}}|
    \leq K_{n} \Big(\sum_{j=1}^{i+1}\Big\{\varphi (y)(y-a)^{2j+1}
+\int_{0}^\infty (v-a)^{2j}\varphi(v)dv\Big\}\Big),\quad 0\leq i
\leq n+1,
\end{equation}
where, from now on, $K_{n}$ denotes  a generic constant, which
only depends on $u$.

\noi Let us introduce
$\widehat{g}_2(\varepsilon):=\widehat{g}(0,0,0,1/\varepsilon), \
\varepsilon \in ]0,1/(2u+2)]$, then :
\begin{equation}\label{asd10}
|\frac{\partial^{i}\widehat{g}_2(\varepsilon)}{\partial
\varepsilon^{i}}|
    \leq K_{n} \int_{0}^\infty v^{2i+2}\varphi(v)dv, \quad 0\leq i
\leq n+1.
\end{equation}
Note that $\varepsilon \leq 1/2$, consequently :
\begin{equation}\label{asd11}
\widehat{g}_2(\varepsilon)\geq \int_0^\infty e^{-v^2/4} \varphi
(v)dv >0.
\end{equation}
Finally, taking into account (\ref{asd9}), (\ref{asd10}),
(\ref{asd11}) and (\ref{int32}) we get :
\begin{eqnarray}
\dis  \Big|\frac{\partial^{n+1}}{\partial \varepsilon ^{n+1}}
\Big(\frac{g(0, a,y,1/\varepsilon
-u)}{g(0,0,0,1/\varepsilon)}\Big)\Big| & \leq  & \dis
K_n\Big(\sum_{j=1}^{n+1}\Big\{\varphi (y)(y-a)^{2j+1}
+\int_{0}^\infty (v-a)^{2j}\varphi(v)dv\Big\}\Big) \nonumber \\
&& \nonumber \\
  |R(0,a,y,u,t)| & \leq  & K_n\dis \Big(\sum_{j=1}^{n+1}\Big\{\varphi (y)(y-a)^{2j+1}
+\int_{0}^\infty (v-a)^{2j}\varphi(v)dv\Big\}\Big) \label{asd12}
\end{eqnarray}

\noi d) Using (\ref{asd1}) and (\ref{asd7}), we have :
$$
\begin{array}{ccl}
  \dis \frac{E_0\big[1_{\Gamma_u}\varphi (S_t)\big]}{E_0[\varphi
(S_t)]} & = & \dis E_0\Big[1_{\Gamma_u}\frac{g(0,X_u,S_u,t-u)}{g(0,0,0,t)}]\\
&&\\
  & = & \dis E_0\Big[1_{\Gamma_u}\Big\{ \sum_{i=0}^n\frac{1}{t^i}F_i(X_u,S_u,u)
  +\frac{1}{t^{n+1}}R(0,X_u,S_u,u,t)\Big\}\Big].\\
\end{array}
$$
Inequality (\ref{asd12}) implies (\ref{int35}) and point 1. b) of
Theorem \ref{tint7}.

\noi e) It remains to prove that for any $i \in \{1, \cdots,n\}$,
$\big(F_i^\varphi(X_t,S_t,t),t\geq 0\big)$ is a $\big( ({\cal
F}_t),P_0)$ martingale.

\noi From (\ref{asd7b}), we deduce that
$E\big[|F_i^\varphi(X_t,S_t,t)|\big]<\infty$.

\noi Let $\Gamma _u\in {\cal F}_u$ and $u\leq v$. The asymptotic
development (\ref{int35}) implies that :
$$\lim_{t\rightarrow \infty}t\Big\{\frac{E_0[1_{\Gamma_u}\varphi (S_t)]}{E_0[\varphi (S_t)]}
-Q_{0}^\varphi (\Gamma _u)\Big\}=
E_0\big[1_{\Gamma_u}F_1^\varphi(X_u,S_u,u)\big] .$$
Since $\Gamma _u\in {\cal F}_v$ then :
$$E_0\big[1_{\Gamma_u}F_1^\varphi(X_u,S_u,u)\big]=
E_0\big[1_{\Gamma_u}F_1^\varphi(X_v,S_v,v)\big].$$
Consequently $\big(F_1^\varphi(X_t,S_t,t),t\geq 0\big)$ is a
$\big( ({\cal F}_t),P_0)$ martingale.

\noi Reasoning by induction, we easily prove  that
$\big((F_i^\varphi(X_t,S_t,t),t\geq 0)\big)$ is a $\big( ({\cal
F}_t),P_0)$ martingale, for any $1\leq i \leq n$.

\noi 3) We now suppose $\lambda >0$. We proceed as previously; the
relation (\ref{asd3}) implies   :
%
%
%
%
%
$$
h(\lambda,z,r) = \sqrt{\frac{2}{\pi }}e^{\lambda ^2
r/2}\Big[\lambda
\sqrt{2\pi}+\frac{1}{\sqrt{r}}e^{-(\frac{z-\lambda r}{\sqrt{r}})
^2 /2}-\lambda \Phi_0\big(\frac{z-\lambda r}{\sqrt{r}}\big)\Big],
$$
where $\Phi_0$  denotes the  function : $\dis
\Phi_0(x):=\int_{-\infty}^x e^{-u^2/2}du$.

\noi Suppose $x<0$. Integrating by parts we have :
$$\Phi_0(x)=\int_{-\infty}^x
\frac{1}{u}ue^{-u^2/2}du=-\frac{1}{x}e^{-x^2/2}- \int_{-\infty}^x
\frac{1}{u^2}e^{-u^2/2}du.$$
Reasoning by induction we can easily prove :
$$\Phi_0(x)=e^{-x^2/2}\sum_{i=0}^n (-1)^{i+1}\frac{a_i}{x^{2i+1}}+
(-1)^{n+1}a_{n+1}\int_{-\infty}^x
\frac{1}{u^{2n+2}}e^{-u^2/2}du,$$
with $a_0=1$ and
$$a_i=1\times 3\times \cdots \times
(2i-1)=\frac{(2i)!}{2^ii!}=E_0[X_1^{2i}],\quad \ i\geq 1.$$

\noi This relation implies :
$$\Phi_0(x)=e^{-x^2/2}\Big[\sum_{i=0}^n (-1)^{i+1}\frac{a_i}{x^{2i+1}}+
O\big(\frac{1}{x^{2n+3}}\big)\Big], \quad x\rightarrow -\infty.$$
Consequently :
$$ h(\lambda, z,r)  =   \sqrt{\frac{2}{\pi }}e^{\lambda ^2
r/2}\Big[\lambda \sqrt{2\pi}+e^{-(\frac{z-\lambda r}{\sqrt{r}}) ^2
/2}\ h_1(\lambda,z,r)\Big]$$
where :
$$
h_1(\lambda, z,r):= -\frac{z}{\lambda r^{3/2}(1-\frac{z}{\lambda
r})} + \sum_{i=1}^n (-1)^{i+1}\frac{a_i}{\lambda^{2i}}
\Big(\frac{1}{\sqrt{r}(1-\frac{z}{\lambda r})}\Big)^{2i+1}
+0\big(\frac{1}{r^{n+3/2}}\big), \quad r\rightarrow \infty.$$
Setting $r=t-u$, where $u>0$ is fixed and $t\rightarrow \infty$,
we get :
$$
\begin{array}{ccl}
  h_1(\lambda, z,t-u) & = & \dis -\frac{z}{\lambda t^{3/2}}
\frac{1}{(1-u/t)^{3/2}(1-\frac{z}{\lambda t (1-u/t)})} \\
   &  &  \\
   & + &\dis  \sum_{i=1}^n \frac{(-1)^{i+1}
a_i}{\lambda^{2i}}\frac{1}{t^{i+1/2}}
\Big(\frac{1}{\sqrt{1-u/t}(1-\frac{z}{\lambda
t(1-u/t)})}\Big)^{2i+1} + 0\big(\frac{1}{t^{n+3/2}}\big), \quad
t \rightarrow \infty,\\
\end{array}
$$
This implies :
$$h_1(\lambda, z,t-u)=\frac{1}{ \sqrt{t}}\Big[
\sum_{i=1}^{n}\alpha_i(\lambda, z,u)\frac{1}{t^i} +
0\big(\frac{1}{t^{n+1}}\big)\Big], \quad t \rightarrow \infty,
$$
where, for any $i$,  $(z,u)\mapsto \alpha_i(\lambda, z,u)$ is a
polynomial function with degree less than $n$, with respect to $z$
or $u$. Moreover we have :
$$\alpha_1(\lambda,
z,u)=-\frac{z}{\lambda}+\frac{1}{\lambda^2}.$$
If $r=t-u$ we have :

$$e^{-(\frac{z-\lambda r}{\sqrt{r}}) ^2
/2}=e^{\lambda z+\lambda ^2u/2}e^{-\lambda
^2t/2}e^{\frac{z^2}{2t(1-u/t)}},$$
therefore :
$$h(\lambda, z,t-u)  =   \sqrt{\frac{2}{\pi }}e^{\lambda ^2
(t-u)/2}\Big[\lambda \sqrt{2\pi}+e^{\lambda z+\lambda
^2u/2}e^{-\lambda ^2t/2}\frac{1}{ \sqrt{t}} \Big\{
\sum_{i=1}^{n}\beta_i(\lambda, z,u)\frac{1}{t^i} +
0\big(\frac{1}{t^{n+1}}\big)\Big\}\Big], \quad t \rightarrow
\infty,
$$
where   $(z,u)\mapsto \beta_i(\lambda, z,u)$ is a polynomial
function with at most degree $n$ with respect to $z$ or $u$.

\noi Note that $\dis \beta_1(\lambda, z,u)=\alpha_1(\lambda,
z,u)=-\frac{z}{\lambda}+\frac{1}{\lambda^2}$

 \noi We are able to come back to relation
(\ref{asd1a}) :
$$\begin{array}{ccl}
  g(\lambda ,a,y,t-u) & = &\dis \sqrt{\frac{2}{\pi }}e^{\lambda ^2 (t-u)/2}\Big[\lambda
\sqrt{2\pi}\Big(\psi
(y)\frac{\sinh\big(\lambda(y-a)\big)}{\lambda}+e^{\lambda
a}\int_y^\infty e^{-\lambda x}\psi (x)dx \Big)  \\
&&\\
   & + & \dis  e^{\lambda
^2u/2}e^{-\lambda ^2t/2}\frac{1}{ \sqrt{t}} \sum_{i=1}^{n}
\gamma_i(\lambda, a,y,u)\frac{1}{t^i} +
0\big(\frac{1}{t^{n+1}}\big)\Big], \quad t \rightarrow \infty, \\
\end{array}
$$
where
$$\gamma_i(\lambda, a,y,u)=\psi (y)e^{\lambda
(y-a)}\int_0^{y-a}e^{-\lambda z} \beta_i(\lambda, z,u)dz+
\int_{y}^\infty \psi(x)\beta_i(\lambda, x-a,u)dx.$$
Note that :
\begin{equation}\label{asd12b}
    |\gamma_i(\lambda, a,y,u)|\leq C(1+u^n)\big(1+|a|^n+\psi (y)e^{\lambda
(y-a)}\big), \quad a\in \mathbb{R},y\geq 0, u\geq 0.
\end{equation}

\noi Introducing :
\begin{equation}\label{asd13}
    F_0^{\lambda,\varphi}(a,y,u):=e^{-\lambda
^2u/2}\Big(\psi
(y)\frac{\sinh\big(\lambda(y-a)\big)}{\lambda}+e^{\lambda
a}\int_y^\infty e^{-\lambda x}\psi (x)dx \Big),
\end{equation}
We have :
$$g(\lambda ,a,y,t-u) = \sqrt{\frac{2}{\pi }}e^{\lambda ^2
t/2}\Big[\lambda \sqrt{2\pi}F_0^{\lambda,\varphi}(a,y,u) +
\frac{e^{-\lambda ^2t/2}}{ \sqrt{t}} \sum_{i=1}^{n}
\gamma_i(\lambda, a,y,u)\frac{1}{t^i} +
0\big(\frac{1}{t^{n+1}}\big)\Big], \quad t \rightarrow \infty,.$$

\noi  Since $\dis F_0^{\lambda,\varphi}(0,0,0)=\int_0^\infty
e^{-\lambda x}\psi
 (x)dx=1$, then
$$\begin{array}{ccl}
 \dis  \frac{g(\lambda ,a,y,t-u)}{g(\lambda ,0,0,t)} & = & \dis
\frac{\dis \lambda \sqrt{2\pi}F_0^{\lambda,\varphi}(a,y,u) +
\frac{e^{-\lambda ^2t/2}}{ \sqrt{t}} \sum_{i=1}^{n}
\gamma_i(\lambda, a,y,u)\frac{1}{t^i} +
0\big(\frac{1}{t^{n+1}}\big)} {\dis \lambda \sqrt{2\pi} +
\frac{e^{-\lambda ^2t/2}}{ \sqrt{t}} \sum_{i=1}^{n}
\gamma_i(\lambda, 0,0,0)\frac{1}{t^i} +
0\big(\frac{1}{t^{n+1}}\big)}
\\
&&\\
   & = & \dis
F_0^{\lambda,\varphi}(a,y,u)+
     \frac{e^{-\lambda ^2t/2}}{ \sqrt{t}} \sum_{i=1}^{n}
F_i^{\lambda,\varphi}( a,y,u)\frac{1}{t^i} +
0\big(\frac{1}{t^{n+1}}\big), \quad t \rightarrow \infty.\\
\end{array}
$$
In particular :
$$F_1^{\lambda,\varphi}(a,y,u)=\frac{1}{\lambda\sqrt{2\pi}}\Big(
\gamma_1(\lambda, a,y,u)-\gamma_1(\lambda,
0,0,0)F_0^{\lambda,\varphi}(a,y,u)\Big).$$
It is easy to compute $\gamma_1(\lambda, a,y,u)$. We have :
$$
\begin{array}{ccl}
  \gamma_1(\lambda, a,y,u) & = & \dis \frac{1}{\lambda ^2}\Big(\psi (y)e^{\lambda
(y-a)}\int_0^{y-a}e^{-\lambda z} (1-\lambda z)dz+ \int_{y}^\infty
\psi(x)\big(1 -\lambda(x-y)-\lambda(y-a)\big)dx\Big)\\
&&\\
   & = & \dis \frac{1}{\lambda ^2}\Big(c(\lambda, \varphi)\varphi_1(y)(y-a)
   + \int_{y}^\infty
\psi(x)\big(1 -\lambda(x-y)\big)dx\Big),\\
\end{array}
$$
with  $\dis c(\lambda, \varphi)=\int_{0}^\infty \psi(x)(1 -\lambda
x)dx$ and
$$\varphi_1(y)=\frac{1}{c(\lambda, \varphi)}
   \Big(\psi(y)-\lambda \int_y^\infty \psi (x)dx\Big), \quad y\geq 0.$$
Setting : %
$$\Phi_1(y):=\int_0^y\varphi_1(z)dz,$$
we obtain :
$$\begin{array}{ccl}
  \Phi_1(+\infty)-\Phi_1(y) & = & \dis \frac{1}{c(\lambda, \varphi)}
   \Big(\int_y ^\infty \psi (x)dx-\lambda \int_y ^\infty \psi(x)(x-y)dx \Big)\\
  &&\\
   & = & \dis \frac{1}{c(\lambda, \varphi)}
   \Big(\int_{y}^\infty
\psi(x)\big(1 -\lambda(x-y)\big)dx\Big). \\
\end{array}
$$
In particular :
$$\Phi_1(+\infty)-\Phi_1(0)=\int_0^\infty \varphi_1(z)dz=1.$$
\noi Moreover :
$$F_1^{\lambda,\varphi}(a,y,u)=\frac{c(\lambda, \varphi)}{\lambda^3\sqrt{2\pi}}\Big(
\varphi_1(y)(y-a)
   + \int_{y}^\infty\varphi_1 (z)dz-F_0^{\lambda,\varphi}(a,y,u)\Big).
$$
This proves point 1. (b) of Theorem \ref{tasd1}.

\end{prooff}

%
\section{ Further discussions about Brownian penalizations}\label{dis}
%

\setcounter{equation}{0} \noi As a conclusion to this paper, we
would like to mention that we are presently developing some
further discussions about Brownian penalizations in three papers
in preparation :
\begin{itemize}
    \item in \cite{RVY4}, we study a number of extensions of Pitman's
    theorem, which are closely related with the penalizations
    found in the present paper;
    \item in \cite{RVY5}, we extend most of the results found in the present
    paper when $(X_t)$ is replaced with $(R_t)$, a Bessel process
    with dimension $d<2$, the weight process being a function of
    the local time of $(R_t)$ at level $0$;
    \item in \cite{RVY6}, we study penalization results for $n$-dimensional
    Brownian motion, when the weight process  is
    $\dis \big(\exp-\int_0^t1_C(X_s)ds\big)$, where $C$
    denotes a cone in $\mathbb{R}^d$ with vertex $0$.
\end{itemize}

\pagebreak
%
%


%





\def\refname{References}
\bibliographystyle{plain}
\bibliography{penalizationstroisHA}
\end{document}